\documentclass[a4paper,10pt]{article}

\usepackage[utf8]{inputenc}

\usepackage{authblk}

\usepackage{ifthen,ifpdf}

\ifpdf
  \usepackage{hyperref}
  \hypersetup{colorlinks=false,allcolors=blue}
  \usepackage{hypcap}
  \pdfpagewidth=\paperwidth
  \pdfpageheight=\paperheight
\fi

\usepackage[utf8]{inputenc}
\usepackage[T1]{fontenc}

\usepackage{amsmath, amsthm, amssymb}
\usepackage{mathtools}
\usepackage{enumerate}
\usepackage{afterpage}
\usepackage{placeins}
\usepackage{tabularx}
\usepackage{dsfont}

\usepackage{graphicx}

\usepackage{algorithm}
\usepackage{algpseudocode}
\usepackage{caption}
\usepackage[labelformat=simple]{subcaption}

\usepackage{color}
\usepackage{units}
\usepackage{array, multirow}
\usepackage{booktabs}

\newcolumntype{M}[1]{>{\vspace{3pt}\raggedleft\arraybackslash}m{#1}}

\usepackage{tikz}
\usepackage{pgfplots}
\usepackage{pgfplotstable}
\usepgfplotslibrary{units}
\pgfplotsset{compat=1.9}
\pgfplotstableset{col sep=comma}
\usepackage{bm}

\usepackage[numbers]{natbib}
\usepackage{anysize} 
\marginsize{2.5cm}{2.5cm}{2cm}{2cm}
\usepackage{stackengine}
\usepackage{wrapfig}
\usepackage{enumitem}
\usepackage[normalem]{ulem}
\theoremstyle{plain}

\theoremstyle{remark}

\makeatletter
\newcommand{\citecomment}[2][]{\citen{#2}#1\citevar}
\newcommand{\citeone}[1]{\citecomment{#1}}
\newcommand{\citewo}[2][]{\citecomment[,~#1]{#2}}
\newcommand{\citevar}{\@ifnextchar\bgroup{;~\citeone}{\@ifnextchar[{;~\citewo}{]}}}
\newcommand{\citefirst}{\@ifnextchar\bgroup{\citeone}{\@ifnextchar[{\citewo}{]}}}

\makeatother

\title{An FFT-based solver with general boundary conditions for stationary diffusion problems based on Chebyshev collocation}

\author[1,2]{J. Quecedo}
\author[2,1,*]{J. Segurado}

\affil[1]{IMDEA Materials Institute, C/ Eric Kandel 2, 28906 Getafe, Madrid, Spain}
\affil[2]{Department of Materials Science, Universidad Politecnica de Madrid, C/ Profesor Aranguren 3, 28040 Madrid, Spain}
\affil[*]{correspondence to: \texttt{javier.segurado@upm.es}}

\date{\today}

\begin{document}

\maketitle
\begin{abstract}
An efficient and robust FFT solver is proposed for diffusion-type problems with general Neumann and Dirichlet boundary conditions, based on a Chebyshev collocation framework. The method combines Chebyshev polynomial approximations with FFT-based operators to provide a matrix-free implementation of the discrete differential operator at the Chebyshev-Gauss–Lobatto points. The linear system of equations resulting from the Chebyshev discretization is solved using LGMRES. To overcome convergence problems on fine grids, a hierarchical refinement strategy based on modal prolongation is proposed, enabling the solution of very large 3D problems. The methodology is available for homogeneous and heterogeneous domains as well for linear and nonlinear consitutive equations.

The accuracy of the proposed method is analyzed by solving the Poisson equation in homogeneous 1D and 3D domains with general boundary conditions, using manufactured analytical solutions as a reference. Convergence to the analytical solution is obtained in a few iterations, reaching smaller errors than those obtained using DCT/DST approaches. Discretizations of up to $256^3$ are achieved, thanks to the hierarchical refinement strategy. In the case of heterogeneous domains, the accuracy and efficiency obtained was similar to that of a standard periodic FFT approach. It is found that the order of complexity of the method preserves the FFT scaling as $n\log n$ in all the cases studied.

\end{abstract}

{\noindent\textbf{Keywords:} Non Periodic Boundary Conditions; Chebyshev collocation method; FFT-based method; Sine--cosine transforms; Hierarchical refinement}

\section{Introduction}
Spectral methods are a well-established class of high-order techniques for the numerical resolution of partial differential equations (PDEs), especially useful when the solution is sufficiently smooth and the global approximation is advantageous \cite{Trefethen,boyd2001,canuto2007,SpectralShen}. 
In contrast to low-order finite-difference or finite-element methods, spectral approaches approximate the solution in terms of global basis functions, such as trigonometric or orthogonal polynomials, and can exhibit very rapid convergence as the approximation order increases \cite{boyd2001,canuto2007,SpectralShen}.

Among spectral techniques, Fourier-based methods have been particularly successful because differential operators can be applied efficiently through the discrete Fourier transform and its fast implementation, the Fast Fourier Transform (FFT) leading to solution methods with computational cost $\mathcal{O}(N \log N)$. 
This feature has made FFT-based schemes especially attractive in computational micromechanics and homogenization where large-scale problems on regular grids can be solved without assembling global sparse or dense matrices \cite{moulinec1994,moulinec1998,suquet2001,vondrejc2013}.  
However, the same trigonometric structure that makes these methods efficient also ties them naturally to periodic settings. 
As a result, the treatment of non-periodic boundary conditions typically requires either reformulation on an enlarged periodic domain or the use of specialized constructions, which may reduce flexibility and complicate the accurate imposition of general Dirichlet or Neumann conditions \cite{Trefethen,boyd2001,SpectralShen}.

Several recent contributions have addressed this limitation while retaining the computational advantages of FFT-based solvers. A first group of approaches, known as framing methods, is based on embedding the domain under study in a larger prism surrounded by a frame or buffer with special properties. The target boundary conditions are prescribed considering this frame by including some condition on the solution scheme. This group includes, among other works, the methods of \cite{Gelebart2020Dirichlet} for Dirichlet boundary conditions or \cite{Zecevic2025ExtendedFFTNonPeriodicBC} who consider both types. An alternative approach is based on applying singular body forces on some internal boundary, as first proposed in 
\cite{Sancho2023ImplicitFFTWave} for elastodynamic problems, in \cite{Zarzoso2025FFTChemoMechanicalFracture} for chemomechanics, and in \cite{buryachenko2025new} for mechanical static approaches. 

A third framework is based on using sine–cosine transform formulations (discrete DST and DST algorithms) as an alternative to standard Fourier transforms to prescribe Dirichlet or Neumann BCs in different problems. In \cite{Gelebart2024ConductingBC,Paux2025DCTSineCosineConductivity}, a solver based on DCT/DST was proposed for conductivity problems   with combined non-uniform Neumann, periodic and Dirichlet boundary conditions. \cite{Risthaus2024ThermalBC} also proposed a method for thermal homogenization with Dirichlet, Neumann and periodic conditions within FFT- and tensor-train-based Green's operator methods .
Risthaus and Schneider developed a direct extended these DCT/DST formulations to enforceme Dirichlet boundary conditions in FFT-based micromechanics, avoiding some of the auxiliary constructions commonly used to enforce prescribed displacements \cite{Risthaus2024DirectDirichlet}.  Also in mechanical problems, Paux et. al. proposed an approach which usses DCT/DST in addition to some conditions on the nodes of the boundary, allowing to impose general boundary conditions, including arbitrary tractions on the boundary 
\cite{Paux2025SineCosineElasticity}.

These developments show the growing interest in extending FFT-based solvers beyond purely periodic settings. 
Nevertheless, most of these approaches remain closely connected to trigonometric approximation spaces and often rely on problem-specific ingredients, such as sine--cosine basis selections, auxiliary decompositions, Green's operators adapted to the boundary conditions, buffer regions, or periodic embeddings. 
Although effective in their respective settings, such strategies do not always provide a unified collocation framework for general non-periodic boundary-value problems on bounded domains.

A natural alternative is to replace purely trigonometric approximations by polynomial spectral discretizations on bounded domains. 
In particular, Chebyshev collocation methods retain the global character and high accuracy of spectral approximations while being well suited to non-periodic problems posed on finite intervals and tensor-product domains \cite{Trefethen,boyd2001,SpectralShen}. 
Moreover, Chebyshev-based differentiation can still be implemented efficiently through fast cosine- and sine-transform techniques, making these methods attractive from a computational standpoint. 
At the same time, the imposition of general boundary conditions and the poor conditioning of high-order differentiation operators remain important practical issues in collocation formulations \cite{hoepffner2007}.

Motivated by these considerations, this work proposes a novel, robust, and efficient Chebyshev collocation framework for diffusion-type problems with fully general boundary conditions on bounded domains. The formulation uses fast transform-based operator evaluations to avoid the explicit assembly of dense differentiation matrices and to obtain a matrix-free implementation. 
A key ingredient of our formulation is the development of a hierarchical refinement procedure, based on modal prolongation, which succeeds in overcoming the convergence limitations of classical Chebyshev implementations and improves the practical performance of the iterative solution process across nested discretizations. The main objective is to provide a practical transform-based Chebyshev framework that extends some of the computational advantages of FFT-enabled spectral methods to non-periodic settings.

The methodology is assessed on representative benchmark Poisson problems, including diffusion on a homogeneous domain, its extension to two-phase heterogeneous domains, and nonlinear heat transfer. 
These examples are used to examine the precision of the formulation, the treatment of boundary conditions, and the practical role of hierarchical refinement in the iterative solution process. 
A comparison with sine--cosine-based transform solvers is also included in order to highlight the relative advantages and limitations of the Chebyshev approach for non-periodic problems.

The remainder of the paper is organized as follows. 
Section~2 presents the numerical formulation, including the Chebyshev collocation discretization, matrix-free differential operators, the treatment of boundary conditions, and the hierarchical refinement strategy. 
Section~3 presents the numerical examples. Finally, Section~4 summarizes the main 
conclusions.
\section{Numerical framework}
Spectral methods aim to solve boundary value problems (BVP) by approximating the solution as a finite sum of known global support smooth functions. In this work, a spectral method is proposed for general stationary boundary-value problems based on Chebyshev polynomials and the FFT algorithm. For simplicity, the formulation in one dimension is presented first. 
\subsection*{Spectral collocation in one dimension}
In a one-dimensional setting, the objective is to find the function $u(x)$ in a spatial domain $D$ with boundary $\partial D$, which fulfills
\begin{eqnarray}
    \mathcal{L}(u(x))+\mathcal{N}(u(x)) = f(x) \qquad x\in D
    \label{eq:general}\\
    B(u(x)) = g(x) \qquad x\in \partial D
    \label{eq:boundaryconditionsgeneral}
\end{eqnarray}
Here, $\mathcal{L} (\cdot)$ and $B(\cdot)$ are linear operators that include 
spatial derivatives. $B$ is defined on the boundary $\partial D$ and imposes the boundary conditions defined as $g(x)$.  $\mathcal{N}(\cdot)$ is a nonlinear operator, and $f(x)$ represents the source term.
The function $u(x)$, solution of the problem in Eqs. (\ref{eq:general}-\ref{eq:boundaryconditionsgeneral}), is an element of the Hilbert space $\mathcal{H}$ with inner product $(\cdot,\cdot)$ and norm $ \parallel \cdot  \parallel$. It belongs to a subspace $\mathcal{B}$ of $\mathcal{H}$ consisting of all functions $u \in \mathcal{H}$ satisfying the boundary conditions $B(u(x))=g(x)$ when $x \in \partial D$.
The approximate solution to the general BVP (Eqs. (\ref{eq:general}-\ref{eq:boundaryconditionsgeneral})) in a finite dimensional function space $P_N$ is expressed as 
\begin{equation}
    u(x) \approx u_N(x) = \sum_{k = 0}^{N}  a_k \phi_k(x)
\end{equation}
where the functions $\phi_n(x)$ are the \emph{basis functions}, belonging to the space $P_N\in \mathcal{H}$. The finite number of unknowns are the coefficients $a_n$.
A residual function $R_N$ can be defined by substituting the solution $u_N$ into \eqref{eq:general}, 
\begin{equation}
    R_N(x) := \mathcal{L}(u_N(x))+\mathcal{N}(u_N(x)) - f(x) 
\end{equation}
and the spectral solution can then be obtained by forcing  the projection of the residual to be zero on a set of \emph{test functions} $\psi_j$.
\begin{equation}
    (R_N, \psi_j)_\omega 
    = \int_D R_N(x) \psi_j(x) \omega(x) \mathrm{d}x = 0, \quad 0 \leq j \leq N 
    \label{eq:residual}
\end{equation}
where \( \omega \) is a positive weight function.
The choice of basis and test functions is one of the main features that distinguish
spectral methods from finite-element and finite-difference methods. In the latter two
methods, these functions are local in character with finite regularities. In contrast, spectral methods employ globally supported smooth functions. 
The most commonly used functions are trigonometric functions or orthogonal
polynomials (typically, the eigenfunctions of Sturm-Liouville problems).

The method proposed is based on a collocation approach. The test functions \( \{ \psi_k \} \) are the Lagrange basis polynomials 
  such that \( \psi_k(x_j) = \delta_{jk} \), where \( \{ x_j \} \) is a set of preassigned collocation points. Hence, $P_N$ is the set of all real algebraic polynomials of degree $\leq N$ and to find the solution $u_N \in P_N$ the residual is forced to zero \( R_N(x_j) = 0 \) at every point in a set of Gauss-Lobatto points $\{x_j\}_{j=0}^N$.  The prescribed conditions become
\begin{equation}
    R_N(x_j) = \mathcal{L}(u_N(x_j))+\mathcal{N}(u_N(x_j)) - f(x_j) = 0 \quad 1 \leq j \leq N-1 
    \label{eq:residualN}
\end{equation}
and at the boundaries $\partial D_-=x_0$ and $\partial D_+=x_N$ satisfies the conditions
\begin{equation}
    \mathcal{B} (u_N(x_0)) = g(x_0),\qquad \mathcal{B} (u_N (x_N)) = g(x_N)
    \label{eq:boundaryN}
\end{equation}
The spectral approach to the BVP in Eq. \eqref{eq:boundaryconditionsgeneral} leads to:
\begin{equation}
\begin{aligned}
\mathcal{L}\!\left(\sum_{k=0}^N a_k\phi_k(x_j)\right)
+
\mathcal{N}\!\left(\sum_{k=0}^N a_k\phi_k(x_j)\right)
=
f(x_j),
\qquad 1\le j\le N-1.
\\
\sum_{k=0}^{N} B_- [a_k\phi_k(x_0)]= g_-, \quad \sum_{k=0}^{N} B_+ [a_k\phi_k(x_N)]= g_+.
\label{eq:problem}
\end{aligned}
\end{equation}

The above system contains \(N+1\) equations and \(N+1\) unknowns \(a_k\). The basis functions \(\phi_k(x)\) are chosen as Jacobi polynomials, whose properties are well known \cite{canuto2007}.

Although the approximation \eqref{eq:problem} has been introduced in modal form, the collocation equations are evaluated at the Gauss--Lobatto nodes. Let
$\{x_j\}_{j=0}^N$, be the collocation nodes and let
$U_j = u_N(x_j)$ denote the nodal values of the interpolant. Equivalently,
$u_N$ can be written in the cardinal Lagrange basis associated with these
nodes as
\begin{equation}
    u_N(x)=\sum_{m=0}^{N} \ell_m(x) U_m ,
\qquad
\ell_m(x_j)=\delta_{jm}.
\end{equation}
Therefore, its derivative evaluated at the collocation nodes is
\begin{equation}
    \left.\frac{\partial u_N}{\partial x}\right|_{x=x_j}
=
\sum_{m=0}^{N} \ell_m'(x_j)\,U_m
\equiv
\sum_{m=0}^{N} D_{jm} U_m,
\qquad j=0,\ldots,N
\end{equation}
where \(D_{jm}=\ell_m'(x_j)\) is the nodal Chebyshev--Gauss--Lobatto
differentiation matrix. As the aim of this work is to develop an efficient method, Jacobi polynomials are selected so that the action of the differentiation operator can be computed using an FFT-type algorithm, with a computational cost of $\mathcal{O}(N \log N)$ instead of the usual \(O(N^2)\) required by a standard differentiation matrix approach.
\subsection*{Extension to the three-dimensional setting}
In three dimensions, we consider the stationary BVP posed on the cube
$D=[a,b]^3$, with $\boldsymbol{x}=(\mathbf{x})\in D$,
\begin{equation}
\mathcal{L}(u(\boldsymbol{x}))+\mathcal{N}(u(\boldsymbol{x})) = f(\boldsymbol{x}) \quad \boldsymbol{x}\in D,
\qquad 
\mathcal{B}(u(\boldsymbol{x}))=g(\boldsymbol{x}) \quad \boldsymbol{x}\in \partial D .
\end{equation}
The spectral collocation approximation is built on the space
$\mathbb P_N(D):=P_N([a,b])\otimes P_N([a,b])\otimes P_N([a,b])$:
\begin{equation}
u(\boldsymbol{x})\approx u_N(\boldsymbol{x})
=\sum_{i,j,k=0}^N a_{ijk}\,\phi_i(x)\phi_j(y)\phi_k(z),
\end{equation}
Let $\{x_p\}_{p=0}^N$, $\{y_q\}_{q=0}^N$, $\{z_r\}_{r=0}^N$ be Gauss--Lobatto nodes in each coordinate direction; the 3D collocation grid is the Cartesian product
$\{(x_p,y_q,z_r)\}_{p,q,r=0}^N$. The collocation formulation enforces
\begin{equation}
R_N(x_p,y_q,z_r)=\mathcal{L}(u_N(x_p,y_q,z_r))+\mathcal{N}\!\left(u_N(x_p,y_q,z_r)\right)-f(x_p,y_q,z_r)=0,
\qquad 1\le p,q,r\le N-1,
\end{equation}
while on boundary nodes the corresponding equations are replaced by the discrete boundary conditions
$\mathcal{B}(u_N(\boldsymbol{x}))=g(\boldsymbol{x})$ evaluated on $\partial D$.
Denote by $U_{pqr}=u_N(x_p,y_q,z_r)$ the nodal values and by
$D^{x},D^{y},D^{z}\in\mathbb{R}^{(N+1)\times(N+1)}$ the one-dimensional differentiation matrices on each corresponding axis. Then the partial derivatives are obtained by separable differentiation:
\begin{equation}
\left[\frac{\partial u_N}{\partial x}\right]_{pqr}
=\sum_{m=0}^N D^{(x)}_{pm}\,U_{mqr},\qquad
\left[\frac{\partial u_N}{\partial y}\right]_{pqr}
=\sum_{m=0}^N D^{(y)}_{qm}\,U_{pmr},\qquad
\left[\frac{\partial u_N}{\partial z}\right]_{pqr}
=\sum_{m=0}^N D^{(z)}_{rm}\,U_{pqm}.
\end{equation}
and, upon vectorization $\boldsymbol{u}=\mathrm{vec}(U)$,
\begin{equation}
\frac{\partial \boldsymbol u}{\partial x}\approx (D^{(x)}\otimes I\otimes I)\boldsymbol{u},\quad
\frac{\partial \boldsymbol u}{\partial y}\approx (I\otimes D^{(y)}\otimes I)\boldsymbol{u},\quad
\frac{\partial \boldsymbol u}{\partial z}\approx (I\otimes I\otimes D^{(z)})\boldsymbol{u},
\end{equation}

\subsection{Chebyshev polynomials as basis functions}
Among the most flexible and thorough options for basis functions with global support are the orthogonal Jacobi polynomials, which will be the foundation of this work. These polynomial families form the core of many practical spectral schemes, and the choice between them typically depends on the geometry of the domain and the desired accuracy and properties of the approximation. In particular, a subtype of Jacobi polynomials, the Chebyshev polynomials, will be used.

\begin{figure}[H]
\centering
\includegraphics[width=0.5\linewidth]{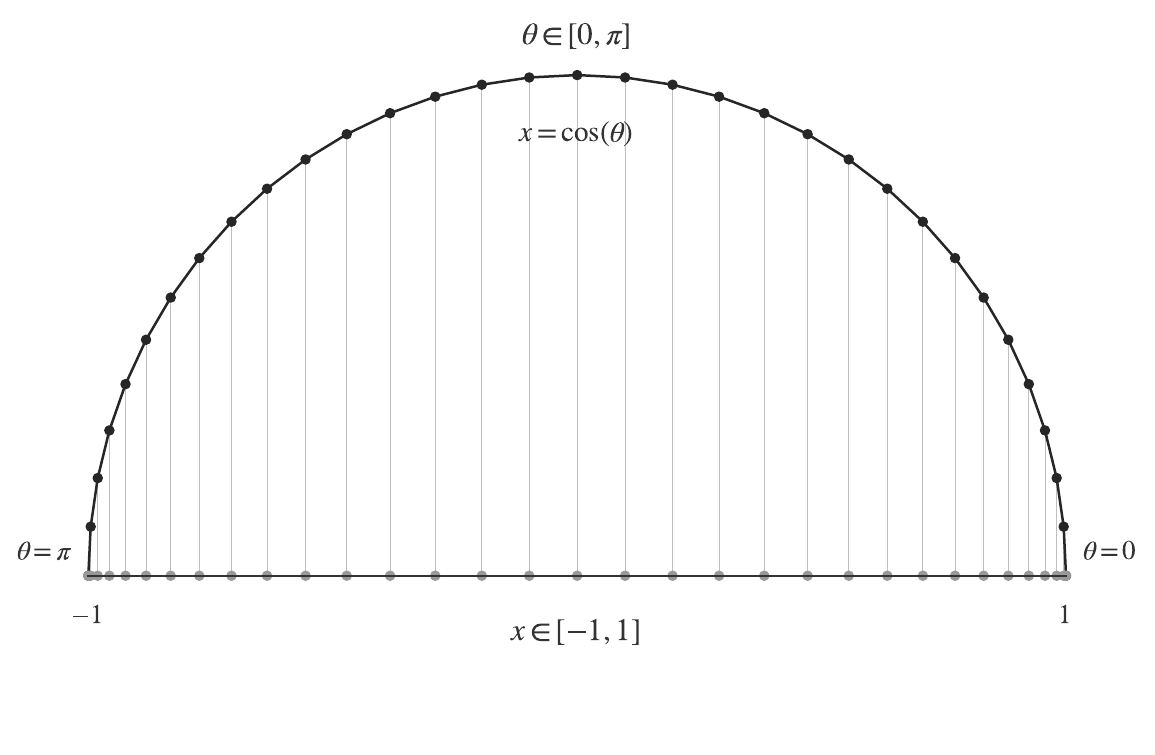}
\captionof{figure}{Projection of uniformly spaced points $\theta$ onto $x$ through $(x=\cos\theta)$, leading to the Chebyshev--Gauss--Lobatto nodes and the cosine representation $(T_n(x)=\cos(n\theta))$.}
\label{fig:circle_complex}
\end{figure}

Chebyshev polynomials admit a simple trigonometric representation, and interpolation and differentiation can be carried out using efficient FFT-based algorithms. This property is highly beneficial as the computational cost tends to $\mathcal{O}(N \log N)$, while being versatile in representing complex solutions. Unlike trigonometric polynomial bases, which are naturally associated with periodic approximations, Chebyshev polynomials provide polynomial approximants on bounded intervals without imposing periodicity or endpoint compatibility conditions. Since algebraic polynomials are dense in spaces of continuous functions on compact intervals, and Chebyshev interpolation exhibits spectral convergence for sufficiently smooth functions, they offer a flexible basis for approximating a wide class of non-periodic solutions while allowing boundary conditions to be imposed independently.
The Chebyshev polynomials are solutions to the differential equation:
\begin{equation}
\sqrt{1 - x^{2}}\;\frac{d}{dx}\!\left(\sqrt{1 - x^{2}}\;\frac{dT_{n}}{dx}\right)
+ n^{2}T_{n}=0,
\qquad x\in[-1,1]
\end{equation}
The analytical expressions of these polynomials are well known and can be obtained from the recurrence relation 
\begin{equation}
\begin{aligned}
& T_0(x) &= 1,\\
& T_1(x) &= x,\\
& T_{n+1}(x)&=2x\,T_{n}(x)-T_{\,n-1}(x).
\end{aligned}
\end{equation}
The Chebyshev polynomials can be expressed in terms of cosines. Let $z$ be a complex number on the unit circle (Fig. \ref{fig:circle_complex}). It can be represented using Euler formula as $z=e^{i\theta} (|z|=1$) where $\theta$ is the argument, and then each $x \in [-1,1]$ can be obtained as
\begin{equation}
    x = Re(z)=\frac{1}{2}(z+z^{-1})=\cos(\theta), \qquad x \in [-1,1].
    \label{eq:x_from_theta}
\end{equation}
$z^{-1}$ is the inverse of $z$ (equal to the conjugate for points in the unit circle). The order $n^{th}$ Chebyshev polynomial can then be expressed in a very simple form as a function of $\theta$ using Eq. \eqref{eq:x_from_theta} as  $T_n(x)=Re(z^n)$
\begin{equation}
\begin{aligned}
 &   T_0(x) = 1 = Re(z^0)=1 \\
 &   T_1(x) = x =  Re(z)=  \frac{1}{2}(e^{j\theta}+e^{-j\theta})=\cos(\theta) \\ 
  & T_2(x) = 2x^2 - 1 = 2(\operatorname{Re}(z))^2-1    = \frac{1}{2}(e^{2j\theta}+e^{-2j\theta}+2)-1=\frac{1}{2}(e^{2j\theta}+e^{-2j\theta})=
   \cos(2\theta) \\
 &  \cdots \\
 &   T_n(x) = \ldots=\cos(n\theta) 
\label{eq:cheb-cosenos}
\end{aligned}
\end{equation}

Equation \eqref{eq:cheb-cosenos} is essential because it allows us to express 
the Chebyshev polynomials in terms of $\cos(n\theta)$, which makes possible to perform calculations using an FFT-based algorithm that implies a cost of order $\mathcal{O}(N \log N)$ \cite{spectral-derivatives-math}.


\subsection{Derivation of the Chebyshev interpolant: FFT-based implementation}
An arbitrary function $f(x)$ in the domain [$-1,1]$ can be approximated by a linear combination of Chebyshev polynomials $T_n$. The approximate expansion of $f(x)$ in the discrete Chebyshev basis, $f_N(x)$, corresponds to
\begin{equation}
    f(x) \approx f_N(x)=\sum_{k=0}^N a_k \space T_k(x) = \sum_{k=0}^N a_k \space \cos(k\theta) =\sum_{k=0}^N a_k \space \cos(k\arccos(x)) 
    \label{eq:f_approx_cos}
\end{equation}
The derivative of $f$ corresponds to
\begin{equation}
\frac{d}{dx}f(\theta(x))
   =
   \frac{d}{d\theta}f(\theta)
   \frac{d\theta}{dx}
   =
   f'(\theta)
   \frac{d}{dx}\!\arccos(x)
   =
   f'(\theta)
   \frac{-1}{\sqrt{1 - x^{2}}}
\end{equation}
with
\begin{equation}
    f'(\theta) = \frac{df(\theta)}{d\theta} \approx \frac{d}{d\theta}\sum_{k=0}^N a_k \space \cos(k\theta)= -\sum_{k=0}^N k a_k \space \sin(k\theta)
    \label{eq:thetaderiv}
\end{equation}
Second order derivatives can be computed in an analogous way
\begin{equation}
    \frac{d^2}{dx^2}f(\theta) = -\frac{d}{dx}(\frac{f'(\theta)}{\sqrt{1-x^2}}) = \frac{f''(\theta)}{1-x^2} - \frac{xf'(\theta)}{(1-x^2)^{3/2}}
    \label{eq:doublethetaderiv}
\end{equation}
Note that the second derivative with respect to $x$ involves both the first and second order derivatives with respect to $\theta$. This effect grows as we move to higher orders: computing the $p$-th derivative with respect to $x$ requires access to every derivative with respect to $\theta$ of order at most $p$.
The term $f'(\theta)$ is calculated as expressed in equation \eqref{eq:thetaderiv} and the term $f''(\theta)$:
\begin{equation}
    f''(\theta) = \frac{d^2f(\theta)}{d\theta^2}= \frac{d^2}{d\theta^2}\sum_{k=0}^N a_k \space \cos(k\theta)= -\sum_{k=0}^N k^2 a_k \space \cos(k\theta)
\end{equation}
Chebyshev polynomials \(T_n(x)\) are defined for \(x \in [-1,1]\). Although differentiating them introduces chain-rule factors involving inverse powers of \(\sqrt{1-x^2}\), the apparent singularities at the endpoints are removable; the cases \(x=\pm 1\) can be treated separately, so the polynomials remain well defined over the entire interval. The expressions of the derivatives at the edges of the domain are well known \cite{Trefethen}. For the first derivative the resulting expressions are
\begin{equation}
    \frac{d}{dx}f(x=1)=\sum_{k=1}^{N} k^2\,a_k
\end{equation}
\begin{equation}
    \frac{d}{dx}f(x=-1)=\sum_{k=1}^{N} (-1)^{k+1} k^2 a_k
\end{equation}
The terms $\frac{d^\mu f(\theta)}{d\theta^\mu}$ can be handled using the Discrete Cosine Transform (DCT) and the Discrete Sine Transform (DST). First, a discrete transformation should be defined to obtain the value of the coefficients $a_k$ as a function of a set of discrete values of the function $f$. 
Considering the uniformly spaced points in $\theta$, $\{x_n=\cos(\theta_n)\}_{n=0}^{N}$, which correspond to the Chebyshev-Gauss-Lobatto points, it is possible to 
obtain the coefficients $a_k$ of the spectral series (Eq. \ref{eq:f_approx_cos}) using the DCT. Among the different implementations of the DCT/DST , the DCT/DST type 1 algorithms from \cite{scipy_dct} are considered. Let be $f_n = f(x_n)$ the values of the function on the Chebyshev-Gauss-Lobatto points, then the coefficients $a_k$ in the cosine expansion resulting in Eq. \ref{eq:cheb-cosenos} are
\begin{equation}
a_k=\frac{c_k}{N}\left[
f_0+(-1)^k f_N+2\sum_{n=1}^{N-1} f_n
\cos\left(\frac{\pi kn}{N}\right)
\right]=\frac{c_k}{N}\textbf{DCT-I}(f),
\qquad k=0,\dots,N,
\end{equation}
and $c_k$ :
\begin{equation}
c_k=
\begin{cases}
\dfrac{1}{2}, & k=0 \text{ or } k=N,\\[4pt]
1, & 1\le k\le N-1.
\end{cases}
\end{equation}

Once the coefficients $a_k$ are obtained, following the expressions \eqref{eq:thetaderiv} and \eqref{eq:doublethetaderiv}, the derivatives of $f(\theta)$ can be calculated using the DST and the DCT, where $\odot$ is the elementwise product.
\begin{equation}
    f'_n = -2\sum_{k=0}^{N}k a_k\sin\left(\frac{\pi k n}{N}\right)=\textbf{IDST-I}(-N(k \odot a_k)), 
    \qquad n = 1,\ldots,N-1
\end{equation}
\begin{equation}
    f''_n = -\big[N^2(-1)^n a_{N}+ \sum_{k=1}^{N-1}k^2a_k \cos\left(\frac{\pi k n }{N}\right)\big]
    = \textbf{IDCT-I}(-N(k^2 \odot a_k)),
    \qquad n = 0,\ldots,N
\end{equation}

The implementation of the latter procedure to compute the $1^{st}$ order derivatives taking advantage of the FFT algorithms is presented 
in Algorithm \ref{alg:ChebDeriv}. 
\begin{algorithm}[H] 
\caption{Chebyshev $1^{st}$ derivative at Chebyshev--Gauss--Lobatto nodes (FFT-based)}
\begin{algorithmic}[1]
\Function{ChebDeriv}{$f,x$} \Comment{$f_n = f(x_n)$ at $x_n=\cos(n\pi/N)$}
    \State $\widehat f \gets \text{DCT-I}(f)$ 
    \Comment{unnormalized type-I discrete cosine transform}
    
    \State $k \gets [0,1,\dots,N-1,0]$
 
    \State $f' \gets 
    \text{IDST-I}\big(-\,(k \odot \widehat f)_{1:N-1}\big)$
    \Comment{$\odot$ elementwise product and inverse type-I sine transform on interior indices}

    \For{$j \gets 1$ to $N-1$}
        \State $(\frac{df}{dx})_j \gets -\dfrac{f'_j}{\sqrt{1-x_j^2}}$
    \EndFor

    \State $df_0 \gets 
    \dfrac{1}{N}\sum_{i=0}^{N-1} i^2\,\widehat f_i 
    \;+\; 
    \dfrac{N}{2}\,\widehat f_N$

    \State $df_N \gets 
    \dfrac{1}{N}\sum_{i=0}^{N-1} (-1)^{i+1} i^2\,\widehat f_i 
    \;+\; 
    \dfrac{N}{2}(-1)^{N+1}\widehat f_N$

    \State \Return $df$
\EndFunction
\end{algorithmic}
\label{alg:ChebDeriv}
\end{algorithm}


Finally, although the Chebyshev construction has been presented on the reference interval $[-1,1]$, the same differentiation procedure can be applied on an arbitrary physical interval $[a,b]$ by means of the affine mapping
\begin{equation}
    t = \frac{b-a}{2}x + \frac{a+b}{2}, 
    \qquad x \in [-1,1], \quad t \in [a,b].
\end{equation}
Therefore, the function must be sampled at the mapped Chebyshev--Gauss--Lobatto nodes
\begin{equation}
    t_n = \frac{b-a}{2}\cos\left(\frac{\pi n}{N}\right) + \frac{a+b}{2},
    \qquad n=0,\ldots,N.
\end{equation}
After computing the derivative on the reference coordinate $x$, the derivative with respect to the physical coordinate $t$ is obtained by the corresponding chain-rule rescaling,
\begin{equation}
    \frac{d^p f}{dt^p}
    =
    \left(\frac{2}{b-a}\right)^p
    \frac{d^p f}{dx^p}.
\end{equation}

The numerical derivative can be easily extended to obtain partial derivatives in higher dimensions. The procedure consists of applying the one-dimensional FFT-based differentiation operator along a prescribed axis, with the remaining coordinate held fixed and with the appropriate affine rescaling from the reference interval $[-1,1]$ to a general physical domain. This yields an efficient matrix-free procedure for evaluating first-order partial derivatives with respect to any selected coordinate direction.

\subsection{Krylov-based PDE Solver}

To solve a PDE with a spectral collocation method, spatial derivatives are replaced by discrete differentiation operators acting on the vector of nodal values, as defined in the previous section. In contrast to local finite-difference stencils, the resulting approximation is \emph{global}, each derivative evaluation at a node depends on the function values at all collocation points. In this work, we apply these operators in a matrix-free approach without explicitly assembling or storing full differentiation matrices.

The general PDE considered is defined in Eq. \eqref{eq:general}. In one spatial dimension, $\mathcal{L}(\cdot)$ may be written as
\begin{equation}
    \mathcal{L}(\cdot) = \sum_{j=0}^{p} l_j \frac{d^j(\cdot)}{dx^j}
\end{equation}
with coefficients $l_j$ and $p$ the maximum order of the spatial derivatives. Let $\{x_i\}_{i=0}^{N}$ be a set of collocation points, and let $\boldsymbol{u} \in\mathbb{R}^{N+1}$ denote the vector of nodal values, $u_i := u(x_i)$. Denote by $\mathcal{C}$ the discrete first-derivative operator acting on nodal vectors, such that $\mathcal{C}(\boldsymbol{u})\approx \partial_x u$ is evaluated at the collocation points. Higher derivatives are obtained by repeated application, $\mathcal{C}^j(\boldsymbol{u})\approx \partial_x^{\,j} u$, or, more generally, by the corresponding discrete operators of order $j$. A discrete version of the operator $\mathcal{L}(\cdot)$ can be defined, $\mathcal{L}_N(\cdot)$, which operates on the vector $\mathbf{u}$ as
\begin{equation}
      \mathcal{L}(u_N)|_{x=x_i} \approx 
      [\mathcal{L}_N(\mathbf{u})]_i := \left[ \sum_{j=0}^{p} l_j\,\mathcal{C}^{j}(\boldsymbol{u})\right]_i
\end{equation}
and \eqref{eq:general} is reduced to a linear system of equations
\begin{equation}
    \left [ \sum_{j=0}^{k} l_j\,\mathcal{C}^{j}(\boldsymbol{u}) \right ]_i \;+\; [\mathcal{N}_N(\mathbf{u})]_i \;=\; f_i, \qquad i = 1,\ldots,N-1
    \label{eq:solveKrylov}
\end{equation}
where $\mathcal{N}_N$ is the discrete version of the nonlinear operator in Eq. \eqref{eq:general}, and $f_i$ denotes the source term sampled at the collocation point $x_i$. 
The extension to multiple spatial dimensions is straightforward by introducing discrete derivative operators in each coordinate direction. For $f:\Omega\subset\mathbb{R}^3\to\mathbb{R}$ sampled on a tensor-product grid, we define operators $\mathcal{C}_x$, $\mathcal{C}_y$, and $\mathcal{C}_z$ such that
\begin{equation}
    \frac{\partial f(\mathbf{x})}{\partial x}\Big|_{\{x_i,y_j,z_k\}}
    \approx \mathcal{C}_x\!\left( \boldsymbol{f} \right)_{ijk},
    \qquad
    \frac{\partial f(\mathbf{x})}{\partial y}\Big|_{\{x_i,y_j,z_k\}}
    \approx \mathcal{C}_y\!\left(\boldsymbol{f}\right)_{ijk},
    \qquad
    \frac{\partial f(\mathbf{x})}{\partial z}\Big|_{\{x_i,y_j,z_k\}}
    \approx \mathcal{C}_z\!\left(\boldsymbol{f}\right)_{ijk},
   \label{eq:partial_discrete} 
\end{equation}
with analogous constructions for higher-order and mixed derivatives.
In the present framework, the discrete residual equation \eqref{eq:solveKrylov} is solved in a matrix-free manner, with all spatial derivatives evaluated through repeated applications of the discrete differentiation operators, thereby avoiding the explicit construction and storage of dense collocation matrices while preserving the global spectral character of the approximation. When \(\mathcal{N} \equiv 0\), the residual is linear and the resulting algebraic system is solved by a Krylov subspace method. When \(\mathcal{N} \neq 0\), the collocation discretization yields a nonlinear algebraic system, and a nonlinear outer iteration is employed, while Krylov methods are used for the linear solves arising at each nonlinear step. In a Picard iteration, the nonlinear term is evaluated at the previous iteration $m$, \(\mathcal{N}(u^{(m)})\), and moved to the right-hand side, so that one solves
\begin{equation}
    \mathcal{L}(u^{(m+1)}) = f - \mathcal{N}(u^{(m)}),
\end{equation}
In a Newton framework, one instead solves a sequence of linearized correction problems
\begin{equation}
J(u^{(m)})\delta u^{(m)}
=
-\left[
\mathcal{L}(u^{(m)})+\mathcal{N}(u^{(m)})-f
\right],
\qquad
u^{(m+1)}=u^{(m)}+\delta u^{(m)}.
\end{equation}
where \(J(\cdot)\) is the Jacobian of the discrete residual.
\subsubsection{Boundary Conditions}

Boundary conditions are enforced as linear constraints on the discrete state, following the procedure of \cite{hoepffner2007}.
Let be the linear system resulting from discretization and collocation (e.g., Eq. \eqref{eq:solveKrylov} in 1D problems),
\begin{equation}
A x=f,\qquad A\in\mathbb{R}^{N+1\times N+1},
\label{eq:linear}
\end{equation}
where $x\in\mathbb{R}^{N+1}$ is the vector of degrees of freedom. Let $n$ be the number of independent boundary conditions, each associated with a given degree of freedom in the boundary $\{x_r^i\}_{i=0}^{n-1}$, and defined as
\begin{equation}
B x=g,\qquad B\in\mathbb{R}^{n\times N+1},\; g\in\mathbb{R}^{n}.
\end{equation}
The unknowns are split into two groups, the ones in which  boundary conditions are applied (removed, $x_r$), and the rest, belonging to the interior of the domain
(interior, $x_k$). Then, the matrix expression of the boundary conditions can be rewritten as
\begin{equation}
\begin{pmatrix} B_k & B_r \end{pmatrix}
\begin{pmatrix} x_k \\ x_r \end{pmatrix}=g
\quad\Rightarrow\quad
x_r = -B_r^{-1}B_k\,x_k + B_r^{-1}g \equiv Gx_k + Hg
\label{eq:BuildSystemBCs}
\end{equation}
with
\(
G=-B_r^{-1}B_k
\)
and
\(
H=B_r^{-1}.
\)
For homogeneous boundary conditions ($g=0$), $x_r=Gx_k$. The practical requirement is that $B_r$ is nonsingular, which is ensured by choosing $x_r$ such that the boundary constraints are linearly independent.

The reordering of unknowns is then applied to the system (Eq. \ref{eq:linear}), leading to
\(
A=\begin{pmatrix}A_{kk}&A_{kr}\\A_{rk}&A_{rr}\end{pmatrix},
\;
f=\begin{pmatrix}f_k\\ f_r\end{pmatrix}.
\)
Substituting $x_r=Gx_k+Hg$ into the system yields a reduced linear system involving only $x_k$,
\begin{equation}
\tilde {A}\,x_k=\tilde f,
\qquad
\tilde A = A_{kk}+A_{kr}G,
\qquad
\tilde f = f_k - A_{kr}Hg.
\label{eq:ReducedSystemBCs}
\end{equation}
After solving for $x_k$, the eliminated variables are recovered by
\begin{equation}
x_r = Gx_k + Hg,
\label{eq:RestrictBCs}
\end{equation}
and the full state $x$ is obtained by undoing the reordering. 

In the present approach, the matrix of coefficients \(A\) is not explicitly assembled, but is provided as a matrix-free differential operator of $\mathcal{L}_N(\cdot)$. 
Although Eqs. \eqref{eq:BuildSystemBCs} and \eqref{eq:ReducedSystemBCs} involve matrix operations, the operations are carried out implicitly through the differential operator, keeping the matrix-free approach
except when it is necessary to calculate the components of $B$ involving Neumann boundary conditions.

\subsubsection{PDE Solver Pseudo-algorithm}
For clarity, the linear PDE solver is presented through a main algorithm describing the global collocation workflow (Algorithm \ref{alg:LinearCollocationSolver}), together with the auxiliary algorithm detailing the matrix-free evaluation of the discrete operator (Algorithm \ref{alg:DiscreteOperator}).
The algorithms given below provide a representative two-dimensional implementation of the spectral collocation solver described in the previous sections. In particular, they illustrate the resolution of the elliptic problem
\begin{equation}
\begin{aligned}
&\frac{\partial}{\partial x}\!\left(-\alpha(x,y)\,\frac{\partial u}{\partial x}\right)
+\frac{\partial}{\partial y}\!\left(-\alpha(x,y)\,\frac{\partial u}{\partial y}\right)
= f(x,y)
\end{aligned}
\end{equation}

\begin{algorithm}[H]
\caption{Linear collocation solver with matrix-free boundary reduction}
\begin{algorithmic}[1]
\Function{LinearCollocationSolver}{$N,\mathrm{bcs},f,L_x,L_y,u_0,\mathrm{rtol},\mathrm{atol},\mathrm{maxiter}$}
    \State $(n_y,n_x) \gets (N)$
    \State $(X,Y) \gets (n_y,n_x,L_x,L_y)$ \Comment{CGL grid}

    \State $f_N \gets f(X,Y)$
    \State $(G,H,g) \gets \text{BuildBoundarySystem}(\mathrm{bcs})$ \Comment{from equation \eqref{eq:BuildSystemBCs}}

    \State $\mathcal{L}_N \gets \textbf{DefineDiscreteOperator}(\mathrm{\cdot})$

    \State $(\widetilde{\mathcal{L}}_N,\widetilde f) \gets$ $\text{EnforceBCs}(\mathcal{L}_N,f_N,G,H,g)$ \Comment{from equation \eqref{eq:ReducedSystemBCs}}

    \State $\widetilde u_N \gets \text{LGMRES}(\widetilde{\mathcal{L}}_N,\widetilde f,u_0,\mathrm{rtol},\mathrm{atol},\mathrm{maxiter})$

    \State $u_N \leftarrow \text{ReconstructFullState}(\tilde{u}_N, G, H, g) $\Comment{from equation \eqref{eq:RestrictBCs}}

    \State \Return $X,Y,u_N$
\EndFunction
\end{algorithmic}
\label{alg:LinearCollocationSolver}
\end{algorithm}

\begin{algorithm}[H]
\caption{Matrix-free evaluation of the discrete linear operator}
\begin{algorithmic}[1]
\Function{DefineDiscreteOperator}{$u$} \Comment{$C_i(\cdot) \leftarrow$ ChebDeriv($\cdot$) to axis=i \text{from Algorithm \ref{alg:ChebDeriv}}}

        \State $u_x \gets \mathcal{C}_x(u)$
        \State $u_y \gets \mathcal{C}_y(u)$
        \State $\mathcal{L}(u) =\Delta u \gets \mathcal{C}_x(-\alpha(x,y) \odot u_x) + \mathcal{C}_y(-\alpha(x,y) \odot u_y)$
            \State \Return $\mathcal{L}(u)$
\EndFunction
\end{algorithmic}
\label{alg:DiscreteOperator}
\end{algorithm}

\subsubsection{A method for convergence improvement: Hierarchical refinement}

In the present collocation framework, the discrete solution is sought in a
finite-dimensional polynomial space built from Chebyshev
polynomials on a domain $\Omega\subset\mathbb{R}^3$. If $N+1$ denotes the
number of Chebyshev--Gauss--Lobatto (CGL) nodes per coordinate direction,
then the corresponding interpolation space is
\begin{equation}
        \mathbb{P}_N := P_{N}\otimes P_{N}\otimes P_{N}.
\end{equation}
Here, $\mathbb{P}_N$ denotes the space of polynomials of degree at most $N$.
Equivalently, in the Chebyshev basis,
\[
P_{N}
=
\operatorname{span}\{T_0,T_1,\ldots,T_N\}.
\]

The numerical solution $u_N\in \mathbb{P}_N$ is defined by enforcing the governing
PDE at the interior CGL collocation points and imposing the boundary
conditions at boundary nodes. Once the nodal values
$\{u_N(x_i,y_j,z_k)\}_{i,j,k=0}^N$ are obtained, the discrete field is fully
determined, since there exists a unique tensor-product polynomial in $\mathbb{P}_N$
interpolating these values.
Hence, evaluating the solution at additional points does not 
introduce any new discretization error, but merely re-evaluates the same discrete polynomial representation.
This property enables a natural hierarchical refinement without re-solving the PDE. The nodal solution is first transformed into modal Chebyshev coefficients
$\{A_{mnl}\}_{m,n,l=0}^N$, and the resulting representation is embedded into
a richer space with $N_f+1$ nodes per coordinate direction,
\begin{equation}
    \mathbb{P}_{N_f} := P_{N_f}\otimes P_{N_f}\otimes P_{N_f}, \qquad N_f>N,
\end{equation}
by zero-padding the coefficient array:
\[
A^{(f)}_{mnl}
=
\begin{cases}
A_{mnl}, & 0\leq m,n,l\leq N,\\
0, & \text{otherwise}.
\end{cases}
\]
Therefore, the refinement is not a new approximation of the PDE solution but an exact prolongation of the already computed spectral 
approximation into a higher-resolution sampling. This is particularly effective in Chebyshev spaces, whose basis functions form a 
rich and highly expressive approximation family for smooth solutions, with fast modal decay and spectral accuracy.

Finally, since the linear systems arising at each refinement level are solved iteratively, the solution on the refined grid is used as the initial guess \(u^{(0)}\) for the solver. In practice, the interpolated field \(u_{N_f}^{(0)}\), obtained by the DCT-based coefficient prolongation described above, is, usually, already very close to the discrete solution in \(\mathbb P_{N_f}\) when the underlying solution is sufficiently smooth. As a result, only a small number of additional iterations is typically required to satisfy the prescribed stopping criterion. This substantially accelerates the solution process across refinement levels and, in many cases, makes it possible to converge reliably even though Chebyshev collocation operators become increasingly ill-conditioned as \(N_f\) grows, without the need for explicit preconditioning. Moreover, because the prolongation is implemented through forward/inverse DCTs and zero-padding, its computational cost is negligible compared with the cost of the iterative solve, so the overall refinement strategy remains highly efficient.

An additional practical consequence of this hierarchical spectral refinement is that it allows us to assess the convergence of the problem, providing a stopping criterion for the refinement. Once \(N\) is sufficiently large, it may become unnecessary to compute further Chebyshev modes by explicitly resolving the problem on a finer grid. Indeed, when the coarse-grid approximation has already entered the asymptotic spectral-convergence regime, the higher-order coefficients neglected at level \(N\) carry negligible energy, and the prolongated field obtained by modal zero-padding provides, up to the prescribed tolerance, essentially the same discrete solution that would be recovered after a new Krylov solve on the refined grid. In practice, this situation can be identified by monitoring the relative variation of a suitable discrete norm, namely
\begin{equation}
e_{\mathrm{ref}}^{(N_f)}
=
\frac{
\left\|
U_{N_f}
-
U_{N\to N_f}
\right\|_2
}{
\left\|U_{N_f}\right\|_2
},
\end{equation}
where \(U_N\) denotes the discrete solution on the current grid and \(U_{N_{f}}\) its prolongation onto the finer grid. When this criterion is satisfied, the approximation may be regarded as numerically converged with respect to further hierarchical refinement, and additional iterative solves on finer meshes become unnecessary. The refinement step then reduces to a purely spectral interpolation procedure, which preserves the already computed solution at negligible cost while avoiding the increasingly ill-conditioned Krylov solves associated with very fine Chebyshev collocation grids.

The procedure of hierarchical refinement is implemented using Algorithm \ref{alg:HierarchicalRefinement}

\begin{algorithm}[H]
\caption{Hierarchical spectral refinement for the linear collocation solver}
\begin{algorithmic}[1]
\Function{HierarchicalRefinement}{$N_0,N_f,\mathrm{bcs},f,\mathrm{coeffs},L_x,L_y,\mathrm{tol},u_0$}
    \State $u \gets \textbf{LinearCollocationSolver}(N_0,\mathrm{bcs},f,\mathrm{coeffs},L_x,L_y,u_0)$

    \State $u^{(0)}_{N_{f}} \gets \text{ChebyshevProlongation}(u,N_{f})$
    \Comment{DCT to modal space, zero-padding, inverse DCT}

    \State $\varepsilon \gets \dfrac{\left\|u^{(0)}_{N_{f}}-u\right\|}{\left\|u\right\|}$

    \State $u \gets \textbf{LinearCollocationSolver}(N_{\mathrm{next}},\mathrm{bcs},f,\mathrm{coeffs},L_x,L_y,u^{(0)}_{N_{f}})$
    \State \Return $u, \varepsilon$
\EndFunction
\end{algorithmic}
\label{alg:HierarchicalRefinement}
\end{algorithm}

\section{Results}
A set of four problems is solved to assess the accuracy, convergence, and efficiency of the matrix-free FFT-based Chebyshev solver proposed. The algorithms 1-4 are programmed in python using \verb|scipy| FFT algorithm and LGMRES solver. 
 The relative and absolute errors for the LGMRES used in all the calculations are $\mathrm{rtol}=10^{-6}$ and $\mathrm{atol}=10^{-8}$, respectively. All simulations are computed in the same Macbook Pro M2 2022 laptop, and computational times refer to this setting.

\subsection{1D Poisson problem: comparison with sine--cosine-based solvers}
\label{sec:comparison_sine_chebyshev}
In this first section, we will use the method developed to solve a one-dimensional Poisson problem with an analytical solution and Dirichlet boundary conditions. The numerical solution for different discretization levels will be compared with the analytical one and the solution provided by sine–cosine transform (DCT/DST)-based methods.

\subsubsection*{Sine--cosine transform solver.}

The sine--cosine transform framework used in \cite{Gelebart2024ConductingBC,Paux2025DCTSineCosineConductivity,DCTDSTGelebart,Risthaus2024ThermalBC,Risthaus2024DirectDirichlet}, among other works, provides an efficient FFT-based strategy to solve transient and stationary diffusion problems on rectangular domains with non-periodic boundary conditions. The main idea is to choose, in each coordinate direction, a trigonometric basis whose symmetry is compatible with the prescribed boundary condition. Dirichlet conditions are naturally associated with sine expansions, while Neumann conditions are associated with cosine expansions. Mixed boundary conditions are handled by using the corresponding half-range sine or cosine transforms.

In this setting, the boundary conditions are embedded directly into the approximation space. Therefore, for homogeneous boundary data, the transformed unknown automatically satisfies the boundary conditions, and no additional boundary equations are required. For non-homogeneous data, one first subtracts a suitable lifting of the boundary values so that the remaining auxiliary problem has homogeneous boundary conditions and can be solved in the same space.

For example, Dirichlet--Dirichlet, Neumann--Neumann, Dirichlet--Neumann, and Neumann--Dirichlet boundary conditions are associated, respectively, with DST-I, DCT-I, DST-III, and DCT-III transforms. Once the appropriate transform has been selected in each coordinate direction, the Laplacian becomes diagonal in modal space. Thus, the solution procedure consists of three steps: applying a forward sine--cosine transform to the right-hand side, solving independent scalar equations for the modal coefficients, and applying the inverse transform.
When fast sine and cosine transforms are used, the resulting computational cost is $\mathcal{O}(N \log N)$.

\subsubsection*{Accuracy and convergence}

To illustrate the behavior of the Chebyshev solver proposed in a simple setting and compare it with the result provided using DCT-DST methods, we consider the one-dimensional boundary-value problem
\begin{equation}
\frac{d^2u}{dx^2}=e^{4x}, 
\qquad x\in[-1,1],
\label{eq:Problem1dDST}
\end{equation}
with homogeneous Dirichlet boundary conditions
\begin{equation}
u(-1)=0,
\qquad
u(1)=0.
\end{equation}
The exact solution of this problem is
\begin{equation}
u(x)=\frac{\exp(4x)-\sinh(4)x-\cosh(4)}{16}.
\label{eq:ejemplo_sines_exacta}
\end{equation}

For this problem, the Dirichlet conditions at both endpoints lead to a DST-I discretization. If \(x_j\), \(j=1,\ldots,N-2\), denote the interior grid points, the numerical solution is represented as a finite sine expansion,
\begin{equation}
u_N(x_j)
=
\sum_{k=1}^{N-2}
\widehat{u}_k
\sin\left(
\frac{k\pi(x_j+1)}{2}
\right).
\end{equation}
Since
\begin{equation}
\frac{d^2}{dx^2}
\sin\left(
\frac{k\pi(x+1)}{2}
\right)
=
-
\left(
\frac{k\pi}{2}
\right)^2
\sin\left(
\frac{k\pi(x+1)}{2}
\right),
\end{equation}
the modal coefficients are obtained directly from
\begin{equation}
\widehat{u}_k
=
-
\frac{\widehat{f}_k}
{\left(k\pi/2\right)^2},
\qquad
k=1,\ldots,N-2,
\end{equation}
where \(\widehat{f}_k\) are the DST-I coefficients of the right-hand side evaluated at the interior points. The numerical solution is then recovered by applying the inverse DST-I.

\begin{figure}[H]
    \centering
    \includegraphics[width=1\linewidth]{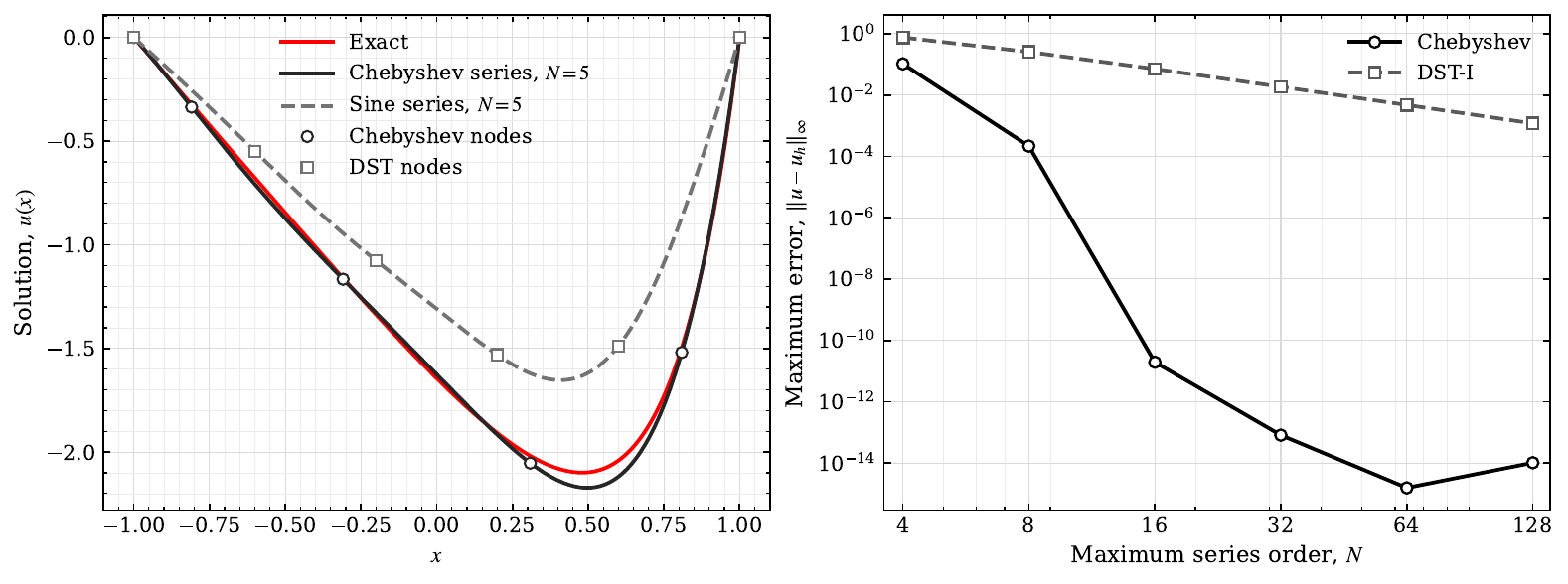}
    \caption{One dimensional Poisson problem: Comparison between the sine-transform approximation and the reference solution.
    (a) Numerical solution for \(N=6\). 
    (b) Maximum error \(\|u-u_h\|_{\infty}\) as a function of the maximum series order.}
    \label{fig:poisson_1d_solution_and_convergence}
\end{figure}

Next, the problem is solved using the Chebyshev method developed, utilizing an increasing discretization ranging from $N=4$ to $N=128$. Fig. ~\ref{fig:poisson_1d_solution_and_convergence}(a) shows the numerical solution for \(N=5\) together with the convergence of the maximum error. The DST-I method enforces the homogeneous Dirichlet boundary conditions exactly, since all sine basis functions vanish at \(x=-1\) and \(x=1\). However, the accuracy is governed not only by the boundary values but also by the smoothness of the corresponding sine expansion.
Chebyshev polynomial expansions constitute a suitable
high-order approximation space for regular solution fields. Their global
polynomial structure, together with their \(C^{N}\)-regularity, allows the
dominant features of smooth solutions to be represented with a relatively
small number of modes. Moreover, once the expansion is truncated, the
resulting approximation space is finite-dimensional; hence, bounded families
in this space exhibit compactness properties under the usual norms are employed
for the solution functions. This provides a stable setting for constructing
convergent modal approximations.

The geometric convergence estimate follows from the classical theory of
Chebyshev approximation for analytic functions; see
~\cite[Theorem~8.2]{trefethen2013approximation}
\begin{equation}
    \|u-u_N^{Ch}\|_{\infty}\leq C\rho^{-N}, \qquad \rho>1.
\end{equation}
By contrast, sine and cosine expansions provide exponential convergence only
when the corresponding odd, even, or periodic extension of the solution remains
sufficiently smooth. If the boundary values or higher derivatives are not
compatible with such an extension, the coefficients decay only
algebraically, leading to estimates of the form
\begin{equation}
    \|u-u_N^S\|\leq C N^{-m}.
\end{equation}
Thus, for non-periodic smooth solutions on bounded intervals, Chebyshev polynomial spaces may represent the solution more efficiently than sine–cosine spaces with the same number of modes, because their convergence is governed by the intrinsic regularity of the solution on the interval rather than by the regularity of an induced periodic extension. To illustrate this convergence behavior, the differences between the results obtained with the Chebyshev and DST methods and the analytical solution are shown in Fig. \ref{fig:poisson_1d_solution_and_convergence}(b). It can be observed that the error decay in the Chebyshev approximation is exponential, reaching machine precision when using around 64 terms. In contrast, the DST approximation decay is algebraic, and for the same refinement level, the error provided by the two approximations differs by 12 orders of magnitude. Regarding the Chebyshev solution, it should be noted that once machine precision is reached, adding more terms becomes useless.

\FloatBarrier
\subsection{3D-Diffusion problem}
\label{subsec:manufactured_linear}

The next case analyzed is a three-dimensional Poisson problem on a unit cube \(\Omega=(0,1)^3\). The objective of this example is to assess
the performance of the matrix-free Chebyshev collocation solver with fine 3D discretization and under mixed Dirichlet--Neumann boundary conditions. In order to evaluate the accuracy of the response, the proposed problem will be obtained by the method of manufactured solutions. The starting point is an analytical expression of a function in the domain,
\begin{equation}
u_{\mathrm{exact}}(\mathbf{x})
=
e^{x+y+z}+x^2+2y-z^2.
\label{eq:exact_solution_3d_example}
\end{equation}
The Laplacian of the select function is \(\Delta u_{\mathrm{exact}}=3e^{x+y+z}\), which will be used as a source term in the Poisson equation. Finally, the boundary conditions on the six faces of the unit cube (named as $\Gamma$ underscored with its corresponding plane) are taken as the local value of the function on the surface points $u_{\mathrm{exact}}$ in the case of Dirichlet, or its gradient projected by the exterior normal in the case of Neumann. The Poisson equation to be solved, including boundary conditions, reads
\begin{subequations}
\label{eq:poisson_3d_example}
\begin{align}
\text{Find}\  \mathbf{u}  \in \Omega \quad \| \quad \Delta u &= 3e^{x+y+z}, \ \text{with} \\
u &= e^{y+z}+2y-z^2,
&& \text{on } \Gamma_{x=0},
\\
u &= e^{1+y+z}+1+2y-z^2,
&& \text{on } \Gamma_{x=1},
\\
\nabla u\cdot n &= -\left(e^{x+z}+2\right),
&& \text{on } \Gamma_{y=0},
\\
\nabla u\cdot n &= e^{x+1+z}+2,
&& \text{on } \Gamma_{y=1},
\\
u &= e^{x+y}+x^2+2y,
&& \text{on } \Gamma_{z=0},
\\
\nabla u\cdot n &= e^{x+y+1}-2,
&& \text{on } \Gamma_{z=1}.
\end{align}
\end{subequations}
In this problem, the function in Eq. \eqref{eq:exact_solution_3d_example} fulfills exactly the differential equation and boundary conditions in Eqs. \eqref{eq:poisson_3d_example}.

The Laplace operator is replaced by the discrete differentiation associated with this problem,
\begin{equation}
\Delta(\cdot)
\approx
\left[
\mathcal{C}_x^2(\cdot)
+
\mathcal{C}_y^2(\cdot)
+
\mathcal{C}_z^2(\cdot)
\right],
\label{eq:poisson_discrete_operator}
\end{equation}
where the directional Chebyshev differentiation operators are defined in Eq. \eqref{eq:partial_discrete} and applied in matrix-free form as
described in Algorithm 1. Boundary conditions are imposed through the elimination procedure introduced in Section~2.4.1, and the resulting reduced linear system is solved with LGMRES.

\begin{figure}[H]
    \centering
    \includegraphics[width=0.8\linewidth]{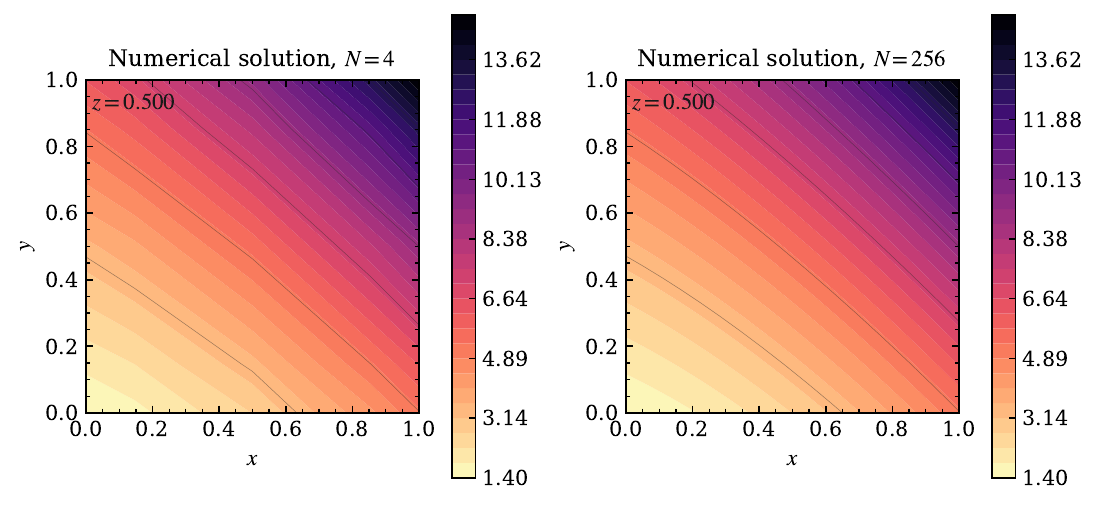}
    \caption{Three dimensional diffusion problem: Numerical solution on the mid-plane \(z=0.5\), comparing (a) low-order and (b) refined Chebyshev collocation approximations.}
\label{fig:poisson_3d_solution_comparison}
    \label{fig:poisson_3d_simple}
\end{figure}

The numerical results obtained in the mid plane are represented in Fig. \ref{fig:poisson_3d_solution_comparison} for two different refinement levels, $N=4$ and $N=256$. It can be shown that, for this relatively smooth problem, the solution obtained by including only 5 terms in the matrix-free Chebyshev approach is qualitatively identical to the one obtained in a highly refined case. 

\begin{figure}[H]
    \centering
    \includegraphics[width=0.85\linewidth]{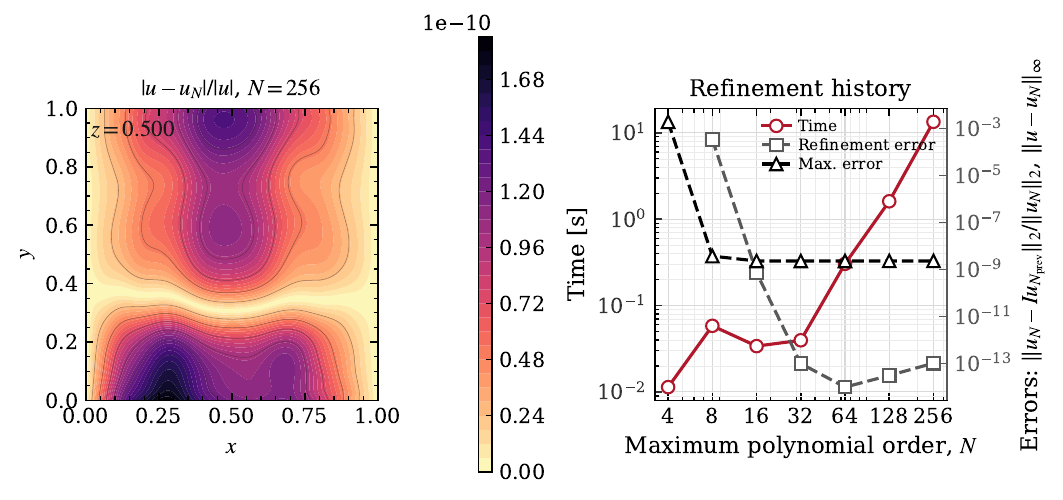}
    \caption{Three dimensional diffusion problem: (a) Error map for $N=256$ on the mid plane with respect analytical solution. (b) Maximum error, refinement error and computational time requiered as function of the Chebyshev polynomial order.}
\label{fig:poisson_3d_error_refinement}
    \label{fig:poisson_3d_refinement}
\end{figure}

To assess the accuracy of the solutions, the local relative error of the numerical solution $u_N$ with respect to the exact one $u$ is obtained as
\begin{equation}
\text{err}(\mathbf{x})=\frac{|u-u_N|}{u}
\label{eq:err_loc}
\end{equation}
and its value is represented in Fig. \ref{fig:poisson_3d_error_refinement}(a) for the same mid-plane, showing values below $10^{-10}$ for $N=256$. This result illustrates the high accuracy of the proposed method. The maximum local error in the full domain with respect to the exact solution is represented in Fig. \ref{fig:poisson_3d_error_refinement}(b), showing its decay to values below $10^{-9}$ for only $N=8$. After that order, the error does not further decrease and is constrained by the tolerance used in the LGMRES solver. To further analyze the effect of the hierarchical refinement, a refinement error is defined as the difference between the prolonged solution and the one obtained in the finer mesh as
\begin{equation}
\text{err}_{ref} = \max \frac{\|u_N-u_{N_{prev}} \|_2}{\|u_N\|_2},
\label{eq:err_ref}
\end{equation}
and its evolution with the polynomial degree is also represented in Fig. \ref{fig:poisson_3d_error_refinement}(b). It can be seen that the refinement error decreases with $N$ (smaller corrections with respect the lower order solution), but when $N=64$ is reached the higher order solution does not change at all with respect to the extrapolated one.
 
Finally, the efficiency of the method is analyzed by representing in Fig. \ref{fig:poisson_3d_error_refinement}(b) the computational time as a function of the number of terms. The tendency observed for $N>32$ is a linear relation between $\log t$ and $\log N$, with slope around 2.7.  If the total number of degrees of freedom of the problem is $n=N^3$, the order of complexity of the FFT would provide $\Delta t\propto n\log n= N^3\log N^3= 3N^3\log N$, and $\log \Delta t \propto 3\log N$, very similar to the slope achieved of 2.7. Therefore, it is quantitatively confirmed that the order of convergence  achieved with the method is $n\log n$, as in any other FFT based approach in homogenization. The number of LGMRES iterations in the hierarchical approach remained quite constant with $N>16$, only requiring one iteration for finer refinements.

\begin{figure}[H]
    \centering
    \includegraphics[width=0.6\linewidth]{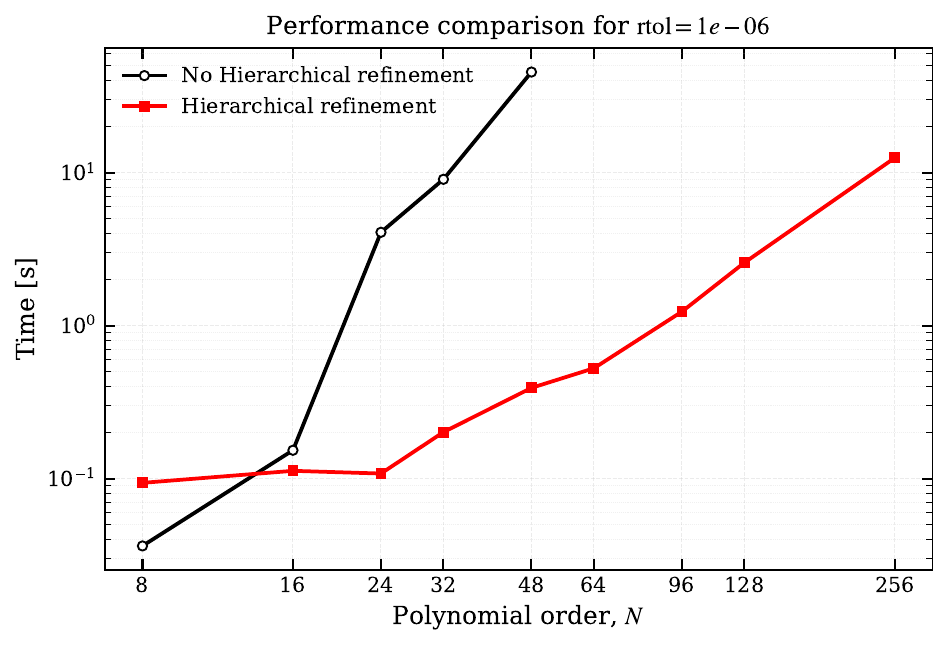}
\caption{Three dimensional diffusion problem: Computational time required as a function of the polynomial order N, with and without hierarchical refinement, using a relative tolerance of $10^{-6}$. }
    \label{fig:comparisonHR}

\end{figure}

The effect of the hierarchical refinement is further analyzed in Fig.   \ref{fig:comparisonHR}, where the computational time is represented as a function of $N$ for the method with and without the use of the hierarchical refinement strategy.

It can be first observed that the rough method fails to converge for $N>48$, due to the increasing ill-conditioning of the Laplace operator with $N$. In contrast, when the hierarchical refinement is introduced, problems with more than 20 million degrees of freedom ($N=256$) are solved without major issues. Moreover, the slope of the computational cost order also decreases with the introduction of hierarchical refinement, indicating a poor scaling of the cost without this technique. The results in Fig. \ref{fig:comparisonHR} highlight the practical relevance of the hierarchical refinement strategy. By using the prolonged solution from a coarser Chebyshev discretization as the initial guess at the next refinement level, the iterative solver starts from a field that already captures the dominant spectral content of the solution. This reduces the computational cost at high polynomial orders, where the conditioning of the collocation operator becomes increasingly demanding while preserving the accuracy of the refined discretization. The comparison therefore shows that hierarchical refinement is not only a post-processing interpolation tool but also an effective acceleration mechanism for matrix-free Chebyshev solvers.


\FloatBarrier
\subsection{3D-Diffusion problem with an spherical inclusion}
Next, a typical homogenization problem is analyzed, since such problems constitute the core application of FFT-based solvers. We consider a diffusion problem in a domain with a position-dependent diffusion coefficient,
\begin{equation}
\nabla\cdot\bigl(-k(\mathbf{x})\nabla u(\mathbf{x})\bigr)=F(\mathbf{x}),
\qquad (\mathbf{x})\in\Omega:=(0,1)^3,
\label{eq:diff_hetero_strong}
\end{equation}
supplemented with boundary conditions that satisfy planar symmetry with respect to the faces of the box,
\begin{align}
& u(0,y,z) = -g \quad \text{on} \ \Gamma_{x=0},\quad \\
& u(1,y,z) = g \quad \text{on} \  \Gamma_{x=1}, \quad \\
& \nabla u \cdot n = 0 \quad \text{on } \ \Gamma \setminus \left(\Gamma_{x=0} \cup \Gamma_{x=1}\right).
\label{eq:bc_mixed_neumann_dirichlet}
\end{align}

Following the matrix-free Chebyshev collocation framework, the directional derivatives are approximated using the operators $\mathcal{C}_x(\cdot)$, $\mathcal{C}_y(\cdot)$, and $\mathcal{C}_z(\cdot)$. At the interior collocation nodes, \eqref{eq1} is discretized as
\begin{equation}
-\left[
\mathcal{C}_x\bigl(k,\mathcal{C}_x(u)\bigr)
+\mathcal{C}_y\bigl(k,\mathcal{C}_y(u)\bigr)
+\mathcal{C}_z\bigl(k,\mathcal{C}_z(u)\bigr)
\right]
= F,
\label{eq1}
\end{equation}

\begin{figure}[H]
    \centering
    \includegraphics[width=0.75\linewidth]{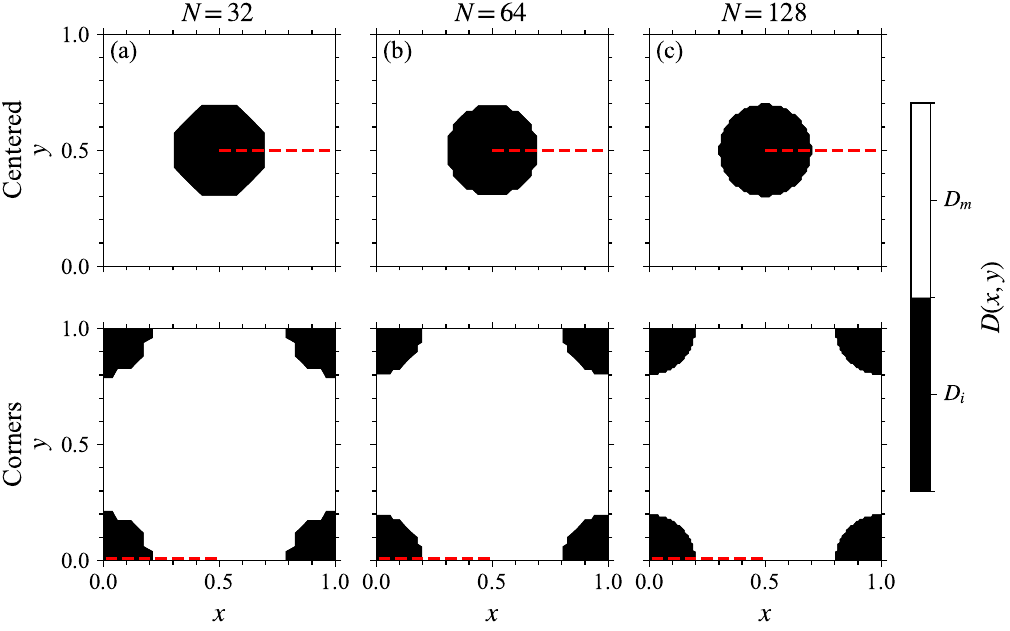}
    \caption{Diffusion with an spherical inclusion: Discretization comparison of the two-phase conductivity field on non-uniform CGL grids.}
    \label{fig:inclusion_center_corner}
\end{figure}

The heterogeneous domain considered here consists of a two-phase composite with a spherical inclusion of radius $R$ and conductivity $k_i$, embedded in a matrix with conductivity $k_m$. The spatially varying diffusion coefficient is defined as
\begin{equation}
k(\mathbf{x})=
\begin{cases}
k_i, & (\mathbf{x})\in\Omega_i,\\
k_m, & (\mathbf{x})\in\Omega\setminus\Omega_i,
\end{cases}
\qquad k_i/k_m = 10,
\label{eq2}
\end{equation}
where $\Omega_i$ denotes the inclusion region. Two inclusion configurations are considered: (i) a single inclusion centered in $\Omega$, and (ii) four quarter-inclusions placed at the corners of $\Omega$, both represented in Fig. \ref{fig:inclusion_center_corner}.
The infinite periodic extension of the domain in the three coordinate directions implies specular symmetry across the six bounding planes. This geometrical symmetry, together with the symmetry of the boundary conditions, makes the solution identical to that obtained under periodic boundary conditions, thereby enabling comparison with standard FFT methods \cite{moulinec1994}.

The two inclusion configurations provide a stringent test of resolution effects intrinsic to Chebyshev--Gauss--Lobatto (CGL) collocation. Although the coefficient field \eqref{eq2} represents the same two-phase medium in both cases, the inclusion representation differs in the two scenarios. CGL nodes are distributed non-uniformly and cluster near the domain boundaries, with the strongest clustering occurring at the grid corners. Consequently, a centered inclusion is predominantly sampled by the relatively coarser interior spacing, whereas corner quarter-inclusions are sampled with a locally finer nodal density. By comparing these two cases, we can assess whether a given CGL resolution provides consistent accuracy for heterogeneities located in the interior and near the corners, and to what extent the non-uniform node distribution affects the representation of sharp coefficient jumps and the resulting solution quality. Nevertheless, note that in a generic microstructure, the local spatial resolution will not be constant as in the case of a regular grid, and the effect of this uneven resolution will depend on the size of the voxels where the interfaces are located. Therefore, in some cases, the local resolution might increase with respect to the regular case, and in others, it might decrease. A conservative metric of the discretization level in a CGL grid with $N+1$ points and length $L=2$ can be the maximum distance between two consecutive points, reached at the center. This value is $h_{CGL}=\cos\left(\frac{\pi}{N}(\frac{N}{2}+1)\right)$, compared to $h_{regular}=\frac{2}{N}$. For a sufficiently large number of points, $h_{CGL}=\frac{\pi}{4}h_{regular} \approx 0.785h_{regular}$, the reduction of the voxel size in the worst position is small.

\begin{figure}[H]
    \centering
    \includegraphics[width=1.0\linewidth]{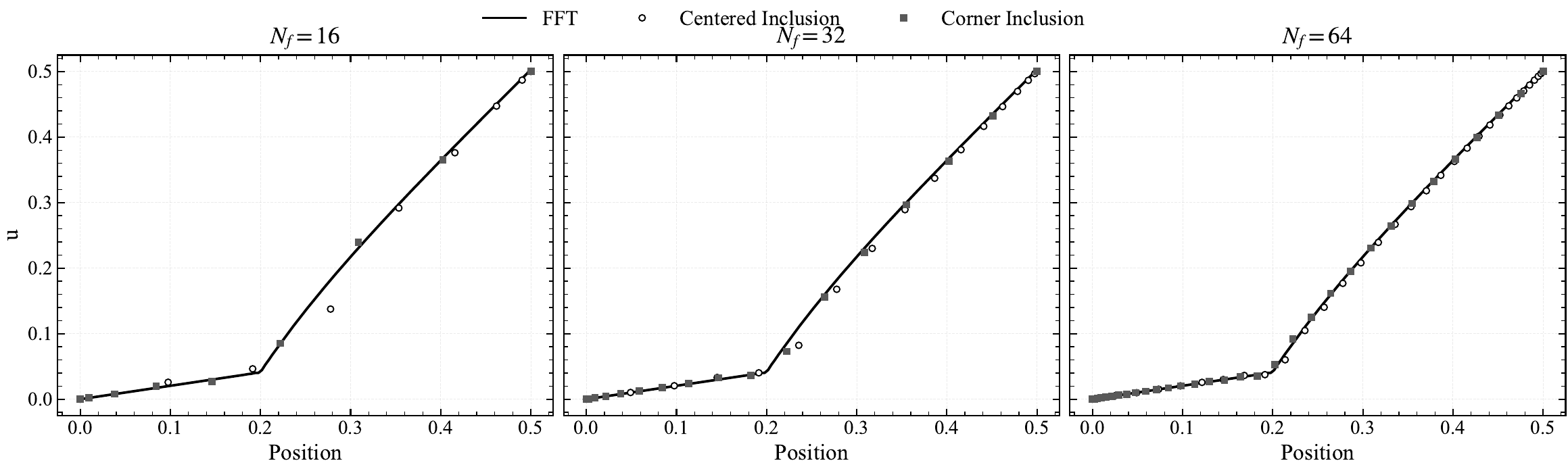}
    \caption{Diffusion with an spherical inclusion: Comparison of $u$  over the dashed red lines of Fig. \ref{fig:inclusion_center_corner} for three final collocation resolutions, $N_f=16$, $32$, and $64$. and corresponding classical FFT solution}
    \label{fig:comparacionFFT}
\end{figure}

The problem \eqref{eq:diff_hetero_strong}--\eqref{eq:bc_mixed_neumann_dirichlet} is solved for different discretization levels, using the Chebyshev method proposed and a standard FFT-based solver (DBFFT approach \cite{LUCARINI2019103131} in the code FFTMAD \cite{lucarini_2018}) on a periodic extension of the unit cell. The value of the field $u$ along the red line depicted in Fig. \ref{fig:inclusion_center_corner}, obtained with the Chebyshev solver and the two particle locations, is represented in Fig. \ref{fig:comparacionFFT} together with the FFT solution for three different discretizations. It can be observed that for very coarse discretization, the location of the sphere influences the local result, and the corner-inclusion case is closer to the FFT profile than the centered-inclusion case, which is consistent with the stronger clustering of CGL nodes near the domain boundaries. However, when the discretization becomes sufficiently fine, the three results become superposed. In this case, using $N=64$ provides an accurate response for the Chebyshev approach. 

\begin{figure}[!htbp]
    \centering
    \includegraphics[width=0.75\linewidth]{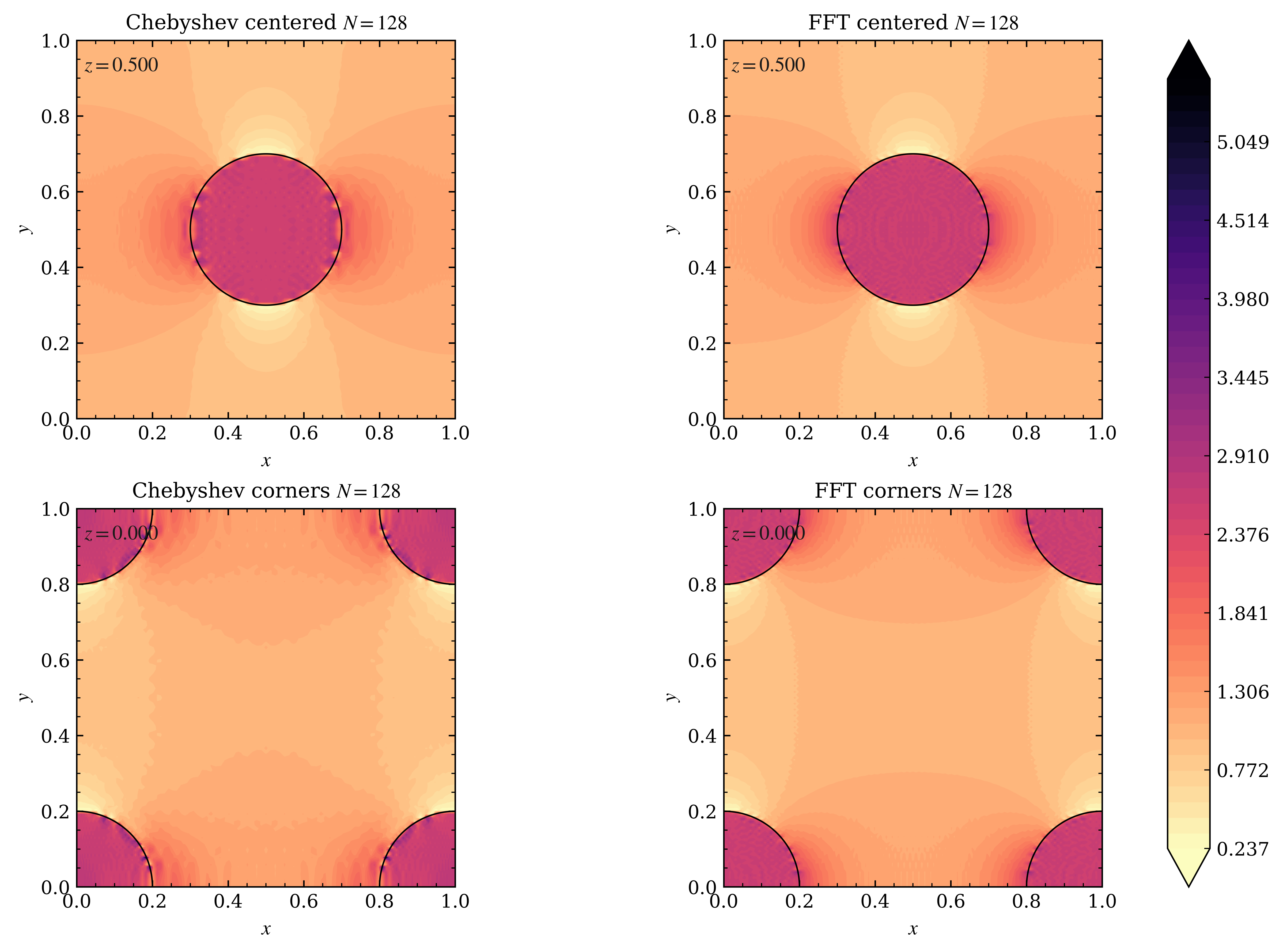}
    \caption{Thermal-flux field for the heterogeneous inclusion problem at (N=128): comparison between the Chebyshev collocation solver and the FFT-based reference solution for centered and corner inclusions.}
    \label{fig:inclusionChebyFFTcenter_corner}
\end{figure}

The resulting thermal-flux in the central plane obtained with the proposed Chebyshev collocation solver and the FFT-based reference solution is represented in Fig.~\ref{fig:inclusionChebyFFTcenter_corner} for the discretization $N=128$. The main differences arise from the distinct discretization strategies used by both methods, which makes the representation of the material interface and the associated numerical oscillations not identical. Despite these local differences, both approaches provide very similar flux distributions, reproducing the same symmetry, localization of the extrema, and overall field pattern around the inclusion. This agreement confirms that the proposed Chebyshev solver captures the dominant response of the heterogeneous diffusion problem while retaining a discretization structure different from that of the FFT reference.

To further assess the robustness of the proposed matrix-free Chebyshev solver on heterogeneous problems, a study of the influence of the phase-property contrast is carried out. In this test, the phase contrast is varied over several orders of magnitude, with \(k_i/k_m \in [10^{-4},10^{4}]\), while keeping the geometrical configuration fixed. For each contrast value, the performance of the LGMRES iterative solver is quantified through the number of iterations required to reach the prescribed tolerance. In addition, the homogenized conductivity obtained from the numerical solution is compared with the  Maxwell  estimation for the dilute limit \cite{Maxwell1873,Kiradjiev2019}.
The effective conductivity predicted by the Maxwell estimate is
\begin{equation}
K_{\mathrm{eff}}^{\mathrm{Maxwell}}
=
K_m
\frac{
2k_m + k_i + 2\phi\left(k_i-k_m\right)
}{
2k_m + k_i - \phi\left(k_i-k_m\right)
}.
\end{equation}
where $\phi$ is the inclusion volume fraction. The numerical effective conductivity is obtained from the volume-averaged heat flux as
\begin{equation}
K_{\mathrm{num}}
=
-\frac{\langle q_y \rangle_{\Omega}}{G},
\qquad
\langle q_y \rangle_{\Omega}
=
\frac{1}{|\Omega|}
\int_{\Omega}
-k(\mathbf{x})\frac{\partial u}{\partial y}\,d\Omega,
\end{equation}
where \(G\) is the imposed macroscopic temperature gradient. For the boundary conditions used in Eq.~(63), this gradient is $G = \frac{2g}{L_y}$.
The relative error in the effective conductivity is then computed as
\begin{equation}
e(K_{\mathrm{eff}})
=
\frac{
\left|K_{\mathrm{num}}-K_{\mathrm{eff}}^{\mathrm{Maxwell}}\right|
}{
\left|K_{\mathrm{eff}}^{\mathrm{Maxwell}}\right|
}.
\end{equation}

\begin{figure}[H]
    \centering
    \includegraphics[width=1\linewidth]{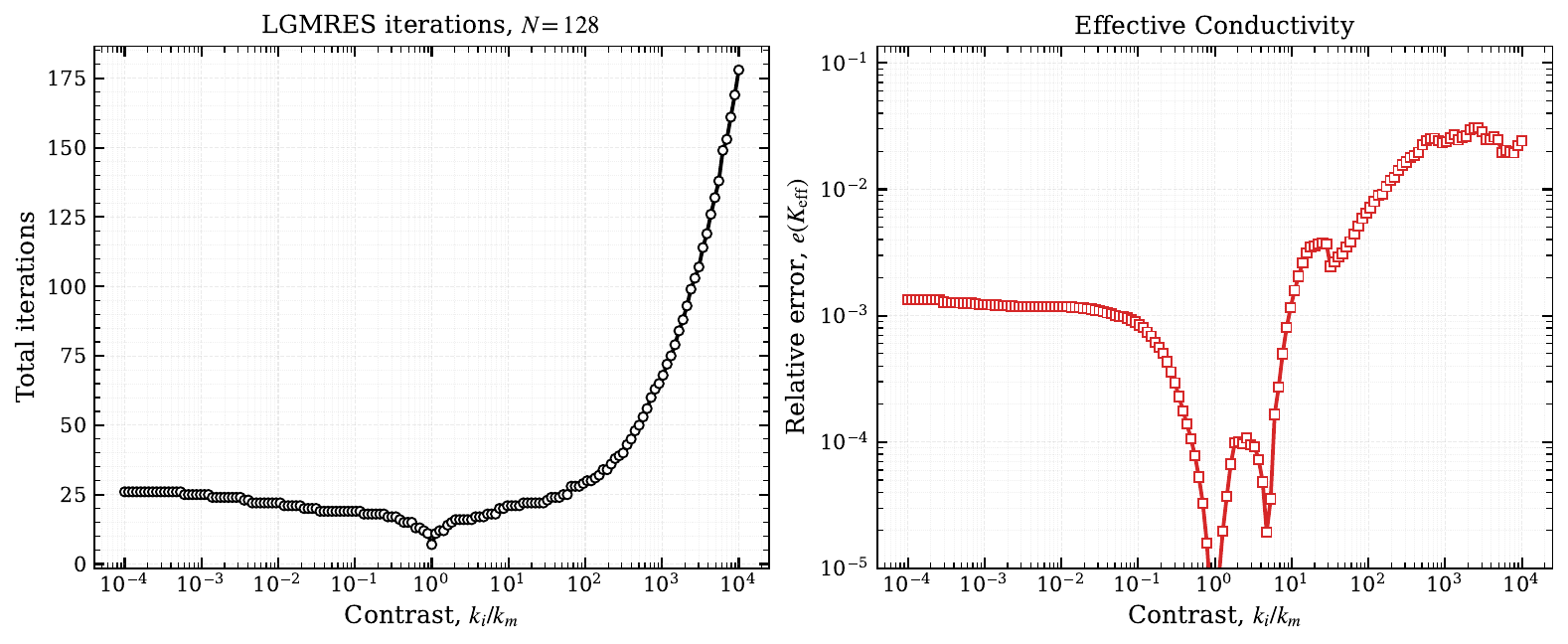}
    \caption{Influence of phase contrast on the iterative convergence and homogenized response of the heterogeneous diffusion problem.}
    \label{fig:contrast_inclusion}
\end{figure}

Fig. \ref{fig:contrast_inclusion}(a) reports the effect of the phase-contrast of the proposed matrix-free Chebyshev solver on its numerical efficiency for a volume fraction $\phi\approx 0.03351$. The results show that the number of LGMRES iterations increases with the phase contrast, as it does in all the different schemes proposed for FFT-based homogenization \cite{lucarini2022}. The convergence deteriorates more when conductivity increases, but the method is able to converge for a wide range of contrasts, making it as versatile as general FFT-based approaches while allowing the use of totally general BCs. Regarding the effective conductivity obtained, the model reproduces the Maxwell estimate with differences below $10^{-3}$ for most of the contrast range. This difference with the Maxwell estimate is reasonable due to the non-dilute character and the regular arrangement of particles in our numerical problem. For very low conductivities, the numerical solution becomes nearly invariant, showing convergence to the case of a hollow sphere. For high conductivities, the solution has slightly more variations but also reaches a relatively constant solution, the infinitely conducting inclusion.

\FloatBarrier
\subsection{3D Nonlinear thermal equilibrium problem}
\label{subsec:manufactured_nonlinear_heat}

As a final demonstration, a nonlinear problem of thermal equilibrium will be solved using the Chebyshev technique proposed. As in previous sections, an exact analytical solution is constructed in order to assess the accuracy and convergence of the numerical solutions obtained with the  method developed.
 In nonlinear heat-conduction problems, the exact temperature field depends explicitly on the constitutive relation adopted for the thermal conductivity; therefore, the manufactured solution must be constructed consistently with the selected nonlinear law.

We consider a nonlinear three-dimensional thermal equilibrium boundary value problem posed on the
unit cube $\Omega=(0,1)^3$ and governed by
\begin{equation}
\label{eq:nl_manufactured_pde}
    \nabla \cdot \left(-k(T)\nabla T\right)=f
    \qquad \text{in } \Omega .
\end{equation}
All quantities in this manufactured problem are written in nondimensional form.
Accordingly, \(T\) denotes a dimensionless temperature, \(k(T)\) denotes the
corresponding normalized thermal conductivity, and \(f\) represents a
dimensionless source term. The thermal conductivity is assumed to depend on the
temperature according to the normalized reciprocal law
\begin{equation}
\label{eq:nl_k_law}
    k(T)=\frac{1}{1+\alpha T},
    \qquad \alpha=0.5 .
\end{equation}
This expression is not intended to represent a universal material law for
metals, but rather a simple admissible nonlinear constitutive relation suitable for verification purposes. Such a monotonic decrease in effective conductivity is consistent with the qualitative behavior
observed in several metallic materials over moderate temperature ranges, where increased scattering at higher temperatures reduces the heat-conduction
capacity. 

To generate a manufactured solution of the nonlinear problem in Eq. \eqref{eq:nl_manufactured_pde}, the first step is to define a sufficiently smooth function $\Phi$. Then, it can be proven that the temperature field
\begin{equation}
\label{eq:exact_solution_3d}
    T(\mathbf{x})
    =
    \frac{\exp\left(\alpha \Phi(\mathbf{x})\right)-1}{\alpha},
\end{equation}
and the source term $f=-\Delta \Phi$ fulfill Eq. \eqref{eq:nl_manufactured_pde}. As in previous examples, the boundary conditions are set based on the values of $T$ and $k(T)\nabla T\cdot \mathbf{n}$ on the surfaces. In this example, we have selected
\begin{equation}
\label{eq:source_3d}
\Phi(\mathbf{x})
    =
    \sin(\pi x)\sin(\pi y)\sin(\pi z) , \qquad
    f
    =
    -\Delta \Phi
    =
    3\pi^2
    \sin(\pi x)\sin(\pi y)\sin(\pi z).
\end{equation}

\begin{figure}[!htbp]
    \centering
    \includegraphics[width=1.0\linewidth]{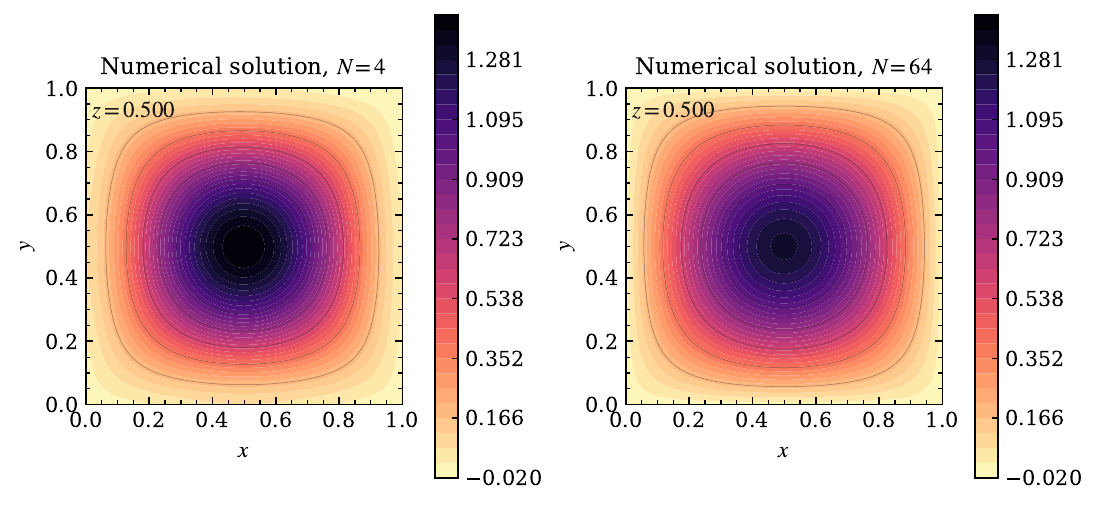}
\caption{Non linear thermal equilibrium: Numerical solution of the nonlinear heat-conduction problem on the mid-plane \(z=0.5\), comparing low-order and refined Chebyshev collocation approximations for \(N=4\) and \(N=128\).}
\label{fig:nonlinear_heat_solution_comparison}
\end{figure}

\begin{figure}[!htbp]
    \centering
    \includegraphics[width=1.0\linewidth]{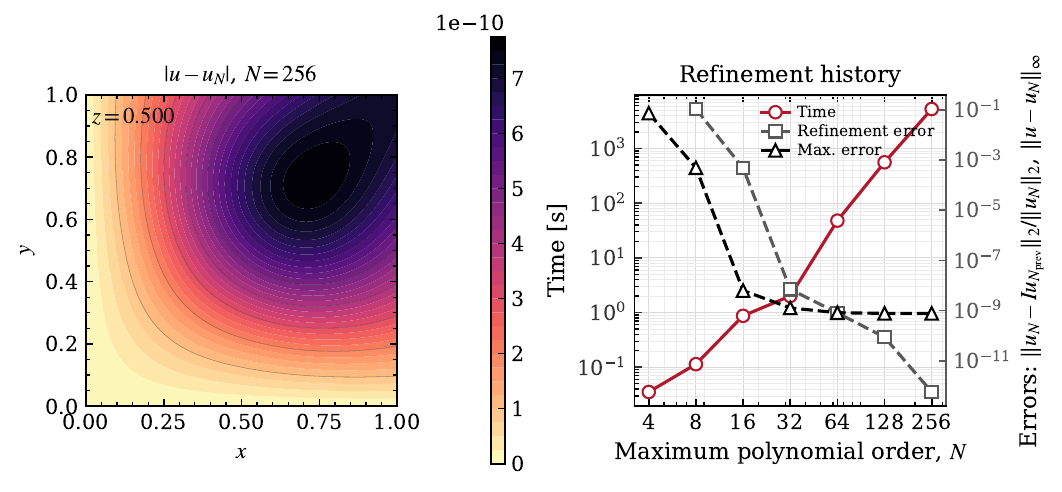}
\caption{Non linear thermal equilibrium: Error distribution and refinement history for the heat-conduction problem, showing the maximum error with respect analytical solution, the hierarchical refinement indicator, and the computational time as the polynomial order increases.}
\label{fig:nonlinear_heat_error_refinement}
\end{figure}

Mixed boundary conditions are imposed consistently with the exact solution. 
Dirichlet data are prescribed on \(\Gamma_{x=0}\) and \(\Gamma_{y=0}\) by
evaluating the manufactured solution on the corresponding boundary faces
\begin{equation}
\label{eq:dirichlet_bcs_3d}
\begin{aligned}
    T(0,y,z)
    &=
    \frac{\exp\left(\alpha \Phi(0,y,z)\right)-1}{\alpha}
    =
    0,
    &&\text{on } \Gamma_{x=0},\\
    T(x,0,z)
    &=
    \frac{\exp\left(\alpha \Phi(x,0,z)\right)-1}{\alpha}
    =
    0,
    &&\text{on } \Gamma_{y=0}.
\end{aligned}
\end{equation}
On the remaining boundary faces, Neumann conditions are imposed in terms of the
nondimensional heat flux,
\begin{equation}
\label{eq:neumann_flux_3d}
    q_N = -k(T)\nabla T\cdot \mathbf n =    -\nabla\Phi\cdot\mathbf n
\end{equation}
where \(\mathbf n\) denotes the outward unit normal. 
Thus, the prescribed Neumann boundary conditions are
\begin{equation}
\label{eq:neumann_data_3d}
\begin{aligned}
    q_N &= -\Phi_x(1,y,z)
    &&\text{on } \Gamma_{x=1},\\
    q_N &= -\Phi_y(x,1,z)
    &&\text{on } \Gamma_{y=1},\\
    q_N &=  \Phi_z(x,y,0)
    &&\text{on } \Gamma_{z=0},\\
    q_N &= -\Phi_z(x,y,1)
    &&\text{on } \Gamma_{z=1}.
\end{aligned}
\end{equation}
With these choices, the temperature field in
Eq.~\eqref{eq:exact_solution_3d} satisfies all the equations of the nonlinear boundary value problem. The nonlinear PDE is solved by Newton's method. Let $R(T)$ be the residual function, 
\begin{equation}
\label{eq:residual_3d}
    R(T)
    =
    \nabla\cdot\left(-k(T)\nabla T\right)-f.
\end{equation} 
The residual is linearized at iteration $m$, corresponding to a field $T^m$, as
\begin{equation}
\label{eq:newton_linearized_problem_3d}
R(T^m+\delta T)\approx R(T^m)+J(T^m)\delta T.
\end{equation}
Forcing the residual to be zero results in the linear equation
\begin{equation}
    J(T^m)\delta T
    =
    -R(T^m), \quad \text{and} \quad  T^{m+1}=T^m+\delta T
    \end{equation}
where the action of the Jacobian over the correction $\delta T$ is
\begin{equation}
\label{eq:jacobian_3d}
    J(T^m)\delta T
    =
    \nabla\cdot
    \left(
        -k(T^m)\nabla \delta T
        -
        k'(T^m)\delta T\nabla T^m
    \right).
\end{equation}
The nonlinear PDE resulting from each Newton increment is solved using the Chebyshev approach, following Algorithm 2 with hierarchical refinement. The dimensionless temperature resulting from two different discretization levels is represented in Fig. \ref{fig:nonlinear_heat_solution_comparison}. 
From a qualitative viewpoint, the results for both both number of polynomial terms and indistinguishable. To assess the validity of the solution, the local error (Eq. \ref{eq:err_loc}) with respect to the analytical solution has been computed for every CGL point and represented in Fig. \ref{fig:nonlinear_heat_solution_comparison} for N=256. It is observed that the error is below $10^{-9}$ at every point again; this accuracy is dictated by the tolerances chosen in the LGMRES. The maximum local error and refinement error (Eq. \ref{eq:err_ref})  are represented on the right side of Fig. \ref{fig:nonlinear_heat_error_refinement}. The conclusions from this graph are the same as for the linear case; the error decays very fast, and a stationary value around $10^{-9}$ is reached from 32 polynomial terms. After that polynomial order, the solution does not improve with respect to the analytical response because the tolerance of the LGMRES is reached. The refinement error, which shows the improvement in the solution by increasing the polynomial order, decays with $N$. The evolution of this metric indicates that the improvement of the solution through the introduction of higher order terms is observed for all entire range of $N$. 

Finally, regarding the computational cost, the right side of Fig. \ref{fig:nonlinear_heat_error_refinement} shows the time spent solving the problem with a different number of terms $N$. The double logarithm diagram has a linear tendency with a slope of 2.98, coinciding almost exactly with the a slope of 3 expected for a scaling depending on the FFT. Therefore, it is quantitatively confirmed that the order of convergence achieved with the method, even in nonlinear cases, is $n\log n$.

To further analyze the numerical performance of the proposed method,  Table~\ref{tbl:nonlinear_hierarchical} presents the number of Newton iterations, LGMRES iterations, errors, and computational time for the different number of polynomial terms $N$.

\begin{table}[H]
\centering
\footnotesize
\setlength{\tabcolsep}{3pt}
\caption{Hierarchical refinement results for the manufactured nonlinear
heat-conduction problem using the matrix-free Chebyshev collocation solver
and Newton--LGMRES iterations. The table reports the Newton iterations, the
total number of LGMRES iterations, the LGMRES iterations per Newton step, the
error with respect to the exact temperature field, the hierarchical refinement
indicator, and the computational time.}
\label{tbl:nonlinear_hierarchical}
\begin{tabular}{@{}r r r l c c r@{}}
\toprule
$N$ &
\shortstack{Newton\\it.} &
\shortstack{LGMRES\\it. total} &
\shortstack{LGMRES it.\\per Newton} &
\shortstack{$\lVert T-T_{\mathrm{exact}}\rVert_{\infty}$} &
\shortstack{$e_{\mathrm{ref}}^{\mathrm{hier}}$} &
\shortstack{Time\\{[s]}} \\
\midrule
4   & 4 & 8  & [2, 2, 2, 2] & $7.352339\cdot10^{-2}$  & --                         & 0.035 \\
8   & 3 & 10 & [4, 3, 3]    & $4.903119\cdot10^{-4}$  & $1.074604\cdot10^{-1}$    & 0.114 \\
16  & 2 & 13 & [7, 6]       & $6.270730\cdot10^{-9}$  & $4.803893\cdot10^{-4}$    & 0.873 \\
32  & 1 & 5  & [5]          & $1.272603\cdot10^{-9}$  & $7.084755\cdot10^{-9}$    & 2.008 \\
64  & 1 & 20 & [20]         & $8.182031\cdot10^{-10}$ & $8.012567\cdot10^{-10}$   & 47.908 \\
128 & 1 & 50 & [50]         & $7.737196\cdot10^{-10}$ & $8.682506\cdot10^{-11}$   & 557.870 \\
256 & 1 & 50 & [50]         & $7.734714\cdot10^{-10}$ & $5.702445\cdot10^{-13}$   & 5245.951 \\
\bottomrule
\end{tabular}
\end{table}

It can be observed that the hierarchical refinement has a dual effect on the performance of the nonlinear problem. On one hand, using a prolonged initial field from an already converged solution at a lower order reduces the number of Newton iterations. Therefore, Newton operates only up to $N=16$, and for higher order, the initial solution is so close to the solution that solving the linearized problem provides the final result. On the other hand, the hierarchical refinement allows the convergence of linear problems with high polynomial orders that would fail to converge without the initial solution. Finally, the time column shows, by taking logarithms, that the $\log \Delta t-\log N$ slope is 3, confirming the FFT convergence order for the whole nonlinear PDE solver.

\FloatBarrier
\section{Conclusions}

In this work, we present a novel approach for extending FFT-enabled spectral solvers to boundary-value problems with fully general, non-periodic boundary conditions, which is based on Chebyshev-based spectral collocation. The proposed method combines Chebyshev polynomial approximations with FFT-based differentiation operators, leading to a matrix-free implementation of the discrete differential operator at the Chebyshev--Gauss--Lobatto points. The resulting linear systems are solved using LGMRES, and to improve robustness and enable the solution of very large problems, a hierarchical refinement strategy based on modal prolongation has also been introduced. The method is also extended to efficiently solve problems with non-linear constitutive equations.

A series of numerical examples is proposed and shows that the method delivers high-fidelity solutions for diffusion, heat-transfer, and related boundary-value problems. In the case of smooth solutions, the method converges to the analytical solution very quickly with the number of polynomial terms. It is shown that the hierarchical refinement strategy constitutes an effective remedy for the ill-conditioning that typically accompanies high-order Chebyshev discretizations. 

The proposed methodology is also able to solve typical problems in homogenization with heterogeneous properties, extended to general boundary conditions and showing similar accuracy and performance to classical FFT homogenization approaches in periodic cases.

Taken together, these results support the proposed approach as a flexible, efficient, and practically scalable spectral methodology with strong potential for the simulation of boundary value problems requiring accurate treatment of non-periodic boundaries.
 
\section*{Acknowledgments}
Javier Segurado acknowledges "Ministerio of Ciencia, Innovación y Universidades" for the project WinPro (Wave Inverse Problems in Heterogeneous media, PID2023-150075OB-I00)

\bibliographystyle{ieeetr}
{\footnotesize
\bibliography{references}
}
\end{document}